\documentclass[english,italian,a4paper]{amsart}
\usepackage          {amssymb} 
\usepackage         {hyperref} 
\usepackage[nobysame]{amsrefs} 
\usepackage            {babel} 
\usepackage   [utf8]{inputenc} 
\usepackage      [T1]{fontenc} 

  \newcommand{\setfont}[1]{\mathbb{#1}} 
  \newcommand{\setC}{\setfont{C}}
  \newcommand{\setF}{\setfont{F}}
  \newcommand{\setQ}{\setfont{Q}}
  \newcommand{\setR}{\setfont{R}}
  \newcommand{\setT}{\setfont{T}}
  \newcommand{\setZ}{\setfont{Z}}
\renewcommand{\leq}{\leqslant}

  \newcommand{\ms}{$'$}

\makeatletter \catcode`\'=11 
\BibSpec{article}{%
    +{}  {\textit}                      {title}
    +{.} { }                            {part}
    +{:} { \textit}                     {subtitle}
    +{,} { \PrintContributions}         {contribution}
    +{.} { \PrintPartials}              {partial}
    +{}  { \PrintAuthors}               {author} 
    +{,} { }                            {journal}
    +{,} { vol.~}                       {volume} 
    +{}  { \PrintDatePV}                {date}
    +{,} { \issuetext}                  {number}
    +{,} { \eprintpages}                {pages}
    +{,} { }                            {status}
    +{,} { \PrintDOI}                   {doi}
    +{,} { available at \eprint}        {eprint}
    +{}  { \parenthesize}               {language}
    +{}  { \PrintTranslation}           {translation}
    +{;} { \PrintReprint}               {reprint}
    +{.} { }                            {note}
    +{.} {}                             {transition}
    +{}  {\SentenceSpace \PrintReviews} {review}
}
\BibSpec{book}{%
    +{}  {\textit}                      {title}
    +{.} { }                            {part}
    +{:} { \textit}                     {subtitle}
    +{}  { \PrintPrimary}               {transition} 
    +{,} { \PrintEdition}               {edition}
    +{}  { \PrintEditorsB}              {editor}
    +{,} { \PrintTranslatorsC}          {translator}
    +{,} { \PrintContributions}         {contribution}
    +{,} { }                            {series}
    +{,} { \voltext}                    {volume}
    +{,} { }                            {publisher}
    +{,} { }                            {organization}
    +{,} { }                            {address}
    +{,} { \PrintDateB}                 {date}
    +{,} { }                            {status}
    +{}  { \parenthesize}               {language}
    +{}  { \PrintTranslation}           {translation}
    +{;} { \PrintReprint}               {reprint}
    +{.} { }                            {note}
    +{.} {}                             {transition}
    +{}  {\SentenceSpace \PrintReviews} {review}
}
\BibSpec{collection.article}{%
    +{}  {\textit}                      {title}
    +{.} { }                            {part}
    +{:} { \textit}                     {subtitle}
    +{}  { \PrintAuthors}               {author} 
    +{,} { \PrintContributions}         {contribution}
    +{,} { \PrintConference}            {conference}
    +{}  {\PrintBook}                   {book}
    +{,} { }                            {booktitle}
    +{,} { \PrintDateB}                 {date}
    +{,} { pagg.~}                      {pages} 
    +{,} { }                            {status}
    +{,} { \PrintDOI}                   {doi}
    +{,} { available at \eprint}        {eprint}
    +{}  { \parenthesize}               {language}
    +{}  { \PrintTranslation}           {translation}
    +{;} { \PrintReprint}               {reprint}
    +{.} { }                            {note}
    +{.} {}                             {transition}
    +{}  {\SentenceSpace \PrintReviews} {review}
}
\renewcommand{\PrintNames@a}[4]{%
    \PrintSeries{\name}
        {(con #1} 
        {}{ e \set@othername} 
        {,}{ \set@othername}
        {}{ e \set@othername} 
        {#2}{#4}{#3)}%
}
\renewcommand{\eprintpages}[1]{%
   pagg.~#1\IfEmptyBibField{eprint}{}{\IfEmptyBibField{journal}{ pp.}{}}
}
\renewcommand{\issuetext}{n.~} 
\makeatother \catcode`\'=12 

\newcommand{\pages}[2]{}
  
\begin{document}

\title
[Una biografia 
di Alessandro Figà-Talamanca (1938--2023)]
{
Una biografia scientifica e personale\\
di Alessandro Figà-Talamanca (1938--2023)\\
e lo sviluppo dell'analisi armonica in Italia}
%
\author
[M.     ~A. Picardello]
{Massimo~A. Picardello}
\address
{Dipartimento di Matematica\\
 Università di Roma ``Tor Vergata''\\
 Via della Ricerca Scientifica\\00133~Roma, Italy}
\email{picard@mat.uniroma2.it}

\maketitle

%


%
%
Alessandro Figà-Talamanca (Roma, 25 maggio 1938 -- Roma, 27 novembre 2023), che gli amici hanno sempre chiamato Sandro (e così faremo qui), è stato un matematico italiano di grande profondità, ma anche di grande capacità di organizzazione della ricerca matematica e della politica scientifica, e soprattutto di comunicazione di idee e di creazione dell'affiatamento ideale affinché in un gruppo di ricerca si creassero forti legami di amicizia. Pertanto queste poche pagine non hanno la pretesa di presentare la profondità della sua opera scientifica, anche se cercherò di tracciare uno schema delle sue idee e di come si siano interconnesse a quelle dei suoi allievi e collaboratori; sono invece un breve cenno della ricchezza della sua personalità e del suo successo nello sviluppare quasi dal nulla in Italia un intero campo di ricerca, quello dell'analisi armonica moderna (in Italia non c'erano mai stati prima studi di analisi armonica, se si eccettua la teoria di Cesàro [1859--1906] della sommabilità di serie, che all'estero aveva avuto applicazioni all'analisi di Fourier, ed una monografia di Tonelli [1885--1946] sulle serie trigonometriche). Quindi, inevitabilmente, questa introduzione è più una breve biografia che un'analisi scientifica, anzi è il racconto di una vita scientifica e di amicizia a lungo in comune con me, e tanti altri, e di come Sandro abbia creato un intero ambito scientifico in Italia, ed allo stesso tempo un grande gruppo di ricerca il quale, grazie alla sua comunicativa ed al suo entusiasmo, fu per molti anni quasi come una famiglia e poi si espanse fino a diventare una scuola. Perciò queste pagine sono un po' anche una biografia di questa ``famiglia'', ma solo per le prime generazioni, altrimenti ci vorrebbe un intero volume; mi scuso del fatto che la presentazione non sia cronologicamente sequenziale bensì un po' intrecciata, nel tentativo di seguire lo sviluppo, in ciascuna sede, di varie generazioni parallele. Mi scuso anche del fatto che la mia collaborazione con Sandro abbia un ruolo preminente in questa biografia: altre collaborazioni sono state per lui ancor più importanti, ma nel delineare sia la sua opera sia la sua energia ed il calore umano della sua personalità sono inevitabilmente guidato dalla sinergia che ho personalmente avuto con lui.

La significativa e continuativa attività nella politica scientifica verrà menzionata solo incidentalmente, perché molti suoi scritti in tale contesto ed in quello parallelo della pubblicistica sono raccolti in altri volumi ed in questa raccolta di scritti matematici non vengono illustrati; ma siccome è indissolubilmente intrecciata all'attività della ricerca scientifica e della sua direzione, ci limitiamo a rammentare che Sandro coprì cariche di grande rilievo nell'organizzazione della matematica: nell'Unione Matematica Italiana (UMI) per oltre un quarantennio, presidente per due mandati (1988--1994) e membro della commissione scientifica per ben undici mandati (1970--1985 e 1994--2012); primo vicepresidente della Società Matematica Europea (1991--1994); anticipo qui che fu fondatore e direttore (1970--1976) della Scuola Matematica Interuniversitaria (SMI), in cui organizzò
i corsi di avviamento alla ricerca tenuti a Perugia e a Cortona da docenti anche stranieri per fornire una prospettiva internazionale ai giovani da avviare alla ricerca, con l'auspicio che al ritorno in Italia percepissero e superassero il diffuso provincialismo del mondo scientifico ed i frequenti egoismi di parte, e mise al contempo in atto un pionieristico sistema di pubblicità capillare dei bandi di concorso per agevolarne il ritorno in Italia; componente del Comitato Nazionale per le Scienze Matematiche del Consiglio Nazionale delle Ricerche (CNR; 1976--1987), del Consiglio Universitario Nazionale (1989--1997) e del Comitato Nazionale per la Valutazione del Sistema Universitario (1999--2004). Fu l'ispiratore del Consorzio Interuniversitario per l’Alta Formazione in Matematica (CIAFM), costituito nel 2004 da vari enti fondatori (Università ed Istituti di ricerca) al fine di finanziare i corsi della SMI, assegnare borse ed assegni di studio per favorire la mobilità dei giovani matematici ed ospitare dottorati in matematica ed altri corsi di studio in forma consortile fra gli enti fondatori. Sandro compì i primi passi presso il Ministero dell'Università per costituire il CIAFM, e Graziano Gentili trovò gli enti disposti a fondarlo (la SNS, la SISSA, l’Università di Perugia, l”INdAM e la SMI); anni dopo, dal 2011 al 2013, Sandro fu il Direttore del CIAFM. Si occupò a lungo anche di altri aspetti politico-organizzativi della ricerca e del sistema universitario anche attraverso un'ampia pubblicistica su alcuni dei maggiori quotidiani e siti nazionali: celebre la sua visionaria e capillare critica all'\emph{Impact Factor}, troppo poco ascoltata, che in un articolo di 33 pagine sul Bollettino dell'UMI nel 1999 previde tutti i danni che la bibliometria ha puntualmente arrecato negli anni e nei decenni successivi. La sua ideazione e realizzazione del trasferimento dei Gruppi Nazionali di ricerca dal CNR all'Istituto Nazionale di Alta Matematica ``F. Severi'' (INdAM) durante la sua presidenza (1995--2003) ha permesso che essi rimanessero gestiti ed amministrati direttamente dalla comunità matematica. Infine, in Sapienza Università di Roma fu direttore del dipartimento di Matematica ``G. Castelnuovo'' (1984--1986 e 2006--2008) e componente della Commissione per l’Innovazione nella Ricerca e Tecnologia (2005--2008), in cui ottenne l'istituzione di percorsi di eccellenza per i migliori studenti.

Questi cenni alla incessante ed infaticabile attività che Sandro ha dedicato al pubblico interesse testimoniano il suo senso dello Stato, e quindi il suo senso etico, che si applica anche alla questione del reclutamento universitario, che ha sempre voluto equo ed aperto anche agli stranieri, e la sua aspirazione alla collaborazione internazionale e l'obiettivo di superare il provincialismo accademico. Vorrei anche aggiungere il suo entusiasmo per l'insegnamento, sia per corsi universitari sia per l'avviamento di allievi alla ricerca, ed il suo costante auspicio che tutti gli allievi nutrissero il suo stesso piacere e la stessa volontà di trasmettere le proprie conoscenze ai propri futuri allievi, e stimolarne la curiosità matematica, anche in campi via via più vari. Sandro manifestava continuamente la propria curiosità per la ricerca matematica e la propria irrefrenabile vivacità, ad esempio tempestando di domande ogni conferenziere. In effetti, ascoltando conferenze o seminari, coglieva o presumeva di cogliere analogie spesso non precise, in qualche modo ispirazioni spontanee, non mediate da una cauta riflessione, collegate ai propri interessi matematici, ed era irrefrenabilmente spinto a porre subito una successione senza tregua di domande, interrompendo spesso il conferenziere. Intendeva infatti chiedere chiarimenti sulla possibilità di trasformare i suoi lampi di intuizione in enunciati precisi da dimostrare rigorosamente, ed auspicava di coinvolgere nel proprio mondo di idee il conferenziere, il quale, poverino, doveva far leva sulla propria maturità scientifica per seguire simultaneamente il filo della propria conferenza ed allo stesso tempo quello, diversissimo, delle idee di Sandro\dots\ Tutti noi che gli eravamo vicini, uditori e conferenziere, sorridevamo e ci divertivamo molto in quei momenti, e non di rado, quando eravamo noi a tenere un seminario od una conferenza, recuperavamo alla fine il tempo perduto continuando l'esposizione in qualche minuto del tempo aggiuntivo dedicato alle domande dell'uditorio, dichiarando che ormai alle domande avevamo già risposto. Naturalmente ci faceva sempre piacere suscitare interesse verso la nostra ricerca e quindi ricevere tante domande grazie alle quali spiegavamo meglio le nostre idee, ma questo era ancora più vero per le domande spontanee ed immediate di Sandro, perché qualche volta l'immediatezza delle sue ispirazioni si rivelava anche foriera di fertili sviluppi. Comunque la vivacità di Sandro non si manifestava solo negli incontri matematici: era ugualmente irrefrenabile anche nei momenti sociali, dove egli manifestava appieno il suo senso dello spirito e la sua vocazione a sviluppare il discorso ed alla narrazione. In effetti, aveva una vasta cultura ed era in grado di affrontare argomenti diversissimi con profonda competenza, ed allo stesso tempo con vivacità. Non si limitava a presentare fatti ed idee già noti, bensì le sue personali rielaborazioni, i suoi punti di vista originali, e spesso nell'ascoltarlo l'interlocutore, anche se di uguale livello di cultura, scopriva approfondimenti nuovi. Sandro aveva un flusso di idee incessante, una chiarezza di idee smagliante, una profonda conoscenza giuridica e politica, una grande capacità di valutare fatti e persone, un forte senso di ciò che è buffo e quindi una naturale tendenza a fare intelligenti battute scherzose, e con slancio irrefrenabile esprimeva tutto questo senza mai alcuna monotonia, in totale sicurezza di sé e con voce stentorea. Non si poteva essere d'accordo con lui su tutto, ma si restava affascinati in ogni caso: le cene con lui erano sempre memorabili.

Sandro era nato a Roma e lì aveva trascorso l'adolescenza, completando la propria educazione secondaria al prestigioso liceo classico Massimo. Cominciò a studiare matematica all'Università di Roma (allora era l'unica università pubblica a Roma, in seguito, quando nel 1979 fu fondata la seconda Università di Roma a Tor Vergata, divenne ``La Sapienza'' e più recentemente, verso fine secolo, ``Sapienza Università di Roma''). Non era molto soddisfatto dei propri studi universitari, ma aveva già sviluppato una notevole curiosità scientifica e vivacità: organizzò un seminario studentesco autogestito nel quale spiegava ai suoi compagni la logica matematica, disciplina a quell'epoca non coltivata a Roma. Per questo fu notato da uno dei suoi professori, Carlo Pucci, il quale, attento e sensibile al potenziale intellettuale dei giovani, lo incoraggiò a continuare gli studi negli Stati Uniti d'America. Fu così che Sandro si trasferì a Los~Angeles e fu ammesso alla Graduate School della University of California at Los~Angeles (UCLA), inizialmente beneficiando di una borsa di studio Fulbright (le borse di studio italiane per studenti all'estero furono varate solo successivamente al CNR da Carlo Pucci con l'entusiastico supporto di Sandro, che grazie ad esse riuscì a far iscrivere intere generazioni di giovani matematici italiani a prestigiose scuole di dottorato estere quando il dottorato di ricerca in Italia ancora non esisteva; sempre fianco a fianco con Pucci, Sandro oltre al programma di borse di studio realizzò, nell'ambito della SMI, un programma di corsi estivi pensato come una leva per indirizzare i giovani brillanti agli studi all'estero). Ad UCLA nel 1964 Sandro conseguì il titolo di Philosophiae Doctor (Ph.D.) sotto la guida di Philip~C. Curtis, Jr.,\ con una tesi, annunciata in~\ocite{01}, sullo spazio $M_p$ dei moltiplicatori dello spazio $L^p(G)$ su un gruppo abeliano $G$ rispetto alla misura di Haar, ovvero delle funzioni che, per moltiplicazione puntuale, preservano lo spazio delle trasformate di Fourier di funzioni in $L^p(G)$: si tratta degli operatori limitati su $L^p(G)$ invarianti per traslazione. Nel caso particolare in cui il gruppo è la retta reale, si tratta dei moltiplicatori di Fourier di $L^p(\setR)$, ossia delle funzioni che, agendo per moltiplicazione puntuale sullo spazio delle trasformate di Fourier di funzioni in $L^p(\setR)$, lo mappano in sé stesso. Più in generale, per ogni $G$ abeliano i casi particolari $p=1$ e $p=2$ erano già noti e più semplici da descrivere. In effetti, $M_1(\Gamma)$ è lo spazio di funzioni definite sul gruppo $\Gamma$ duale di $G$, ovvero il gruppo dei caratteri di $G$, che sono trasformate di Fourier delle misure di Borel regolari e limitate sul gruppo $G$, mentre $M_2(\Gamma)$ è lo spazio $L^\infty(\Gamma)$ delle funzioni limitate e misurabili (nel seguito, rispettando la terminologia tradizionale, invece che $M_1=M_1(\Gamma)$ talvolta scriveremo $B=B(\Gamma)$: è lo spazio delle combinazioni lineari delle funzioni di tipo positivo sul gruppo $\Gamma$). Si noti quindi che $M_1(\Gamma)$ risulta isomorfo allo spazio duale dello spazio $C_0(G)$ delle funzioni continue che tendono a zero all'infinito, ossia al duale dello spazio delle combinazioni lineari di prodotti di convoluzione del tipo $L^1(G)*C_0(G)$ (nella norma indotta naturalmente dai fattori). D'altra parte, $M_2(\Gamma)$ è isomorfo al duale di $L^1(\Gamma)$, cioè al duale dello spazio $A(G)$ generato dalle combinazioni lineari dei prodotti di convoluzione del tipo $L^2(G)*L^2(G)$ (sempre nella norma indotta naturalmente dai fattori). Nella sua tesi Sandro dimostrò che, per $1<p<\infty$, $M_p(\Gamma)$ è il duale dello spazio $A_p(G)$ generato dalle convoluzioni del tipo $L^p(G)*L^q(G)$, con $1/p+1/q=1$, e diede come conseguenza alcune caratterizzazioni di insiemi lacunari su $\Gamma$, estendendo una tematica classica dell'analisi di Fourier su $\setR$ e $\setT$. Tutte le dimostrazioni furono pubblicate in~\ocite{02}. In seguito, Carl Herz dimostrò che gli spazi $A_p(G)$ sono algebre, ossia sono chiusi rispetto alla moltiplicazione puntuale, e da allora furono universalmente noti come \emph{algebre di Figà-Talamanca -- Herz}, (o anche \emph{algebre di Figà-Talamanca}, in Encyclopedia of Mathematics), e divennero un importante ambito di ricerca ed uno strumento indispensabile per lo studio di operatori di convoluzione su $L^p$ di un gruppo localmente compatto.

In questo modo, la linea di ricerca di Sandro si inseriva in un ambito di ricerca in pieno sviluppo, dandovi un forte impulso: quello delle proprietà classiche dell'analisi di Fourier estese a gruppi localmente compatti abeliani. Lo studio dell'analisi su questi gruppi, le cui basi (gruppo duale, formula di Plancherel) erano state introdotte da Pontrjagin (1934) e si erano sviluppate con André Weil (1940), si stava evolvendo nella direzione dell'analisi di Fourier, grazie all'opera di importanti precursori (Walter Rudin, Robert E.~Edwards, Edwin Hewitt ed altri). Il libro \textit{Fourier Analysis on Groups}, di Walter Rudin, era stato pubblicato solo due anni prima della tesi di Sandro.

Più in generale, era stata dedicata attenzione all'analisi armonica su gruppi localmente compatti, non solo commutativi: il libro di Lynn~H. Loomis, \textit{An Introduction to Abstract Harmonic Analysis}, era apparso nel 1953. In effetti, in quel periodo erano di particolare rilevanza per la fisica non i gruppi abeliani ma i gruppi compatti e i gruppi di Lie (ad esempio, per la meccanica quantistica il gruppo di Heisenberg e per la teoria della relatività il gruppo $SO(3,1)$, ed i gruppi $SL(2,\setC)$ e $SL(2,\setR)$ ad esso correlati) e quindi le loro rappresentazioni unitarie. La teoria delle rappresentazioni dei gruppi compatti, con particolare riferimento ai gruppi di Lie, era stata studiata da Hermann Weyl. La teoria delle rappresentazioni dei gruppi di Lie semisemplici e le corrispondenti funzioni sferiche erano l'ambito di studio di un gran numero di famosi matematici, ad esempio a Princeton (la scuola di Harish-Chandra e Robert~P. Langlands) ed a Mosca (la scuola di Izrail\ms~M. Gel\ms{}fand). In ambiti simili a questi Sandro avrebbe in seguito fornito contributi scientifici, ma, in linea con le tematiche sviluppate nella propria tesi, orientò quasi sempre la sua attenzione verso lo studio dei moltiplicatori di Fourier di $L^p$, degli insiemi lacunari e delle funzioni di tipo positivo più che sulla teoria delle rappresentazioni --- con la notevole eccezione, molti anni dopo, rappresentata dalle funzioni sferiche, non su gruppi di Lie bensì su gruppi discreti di grande importanza.

Sandro sviluppò ulteriormente ed espanse grandemente queste linee di ricerca negli anni successivi, in una serie di articoli in collaborazione con Garth Gaudry \cites{06,12,14,16} sui moltiplicatori di $L^p$ su gruppi abeliani e con Daniel Rider \cite{05} sugli insiemi lacunari sul duale di gruppi compatti non necessariamente abeliani. Il duale di un gruppo commutativo (l'insieme dei caratteri, ovvero la famiglia delle sue rappresentazioni unitarie irriducibili a meno di equivalenza), è un gruppo, ma non lo è il duale di un gruppo $G$ non commutativo: quindi, con le ricerche insieme a Rider, Sandro allargava i propri interessi dall'analisi di Fourier sui gruppi all'analisi armonica astratta; inoltre, sempre con Rider, studiò l'appartenenza ad $L^p$ con probabilità 1 di serie di Fourier aleatorie su gruppi compatti \cite{07}. Questa collaborazione fu sviluppata al Massachusetts Institute of Technology (MIT), dove Sandro era stato nominato Moore Instructor per il biennio 1964--1966. Là nel 1965 Sandro conobbe Irene Petritsi, allora ricercatrice alla Harvard School of Public Health, con cui si sposò alla fine di quell'anno e che gli rimase compagna per tutta la vita. Nel biennio 1966--1968 Sandro insegnò all'Università di Genova, e poi nel 1968--1969 fu Lecturer a Berkeley e nel 1969--1970 a Yale, dove Gaudry era Gibbs Instructor. La collaborazione con Gaudry si indirizzò a vari problemi collegati ai moltiplicatori di $L^p$, precisamente sugli insiemi di unicità di $L^p$ in gruppi abeliani e sulla costruzione di moltiplicatori che tendono a zero all'infinito in $\setR^n$ o $\setZ^n$ ma non sono nella chiusura di $A(\setR^n)$ o delle funzioni a supporto finito su $\setZ^n$. Attraverso questa collaborazione Sandro ebbe un influsso, diretto o indiretto, sulla formazione di alcuni giovani matematici australiani. Nel 1970 Sandro lascio l'America e ritornò all'Università di Genova, dove nacque il suo primo figlio Giovanni e dove egli divenne Professore e cominciò ad avviare studenti alla ricerca matematica. Aveva già fin da allora una straordinaria comunicativa con i suoi studenti, ed in seguito vedemmo anche come sapesse adattarla alla psicologia di ciascuno per raggiungere la massima empatia. Il primo di essi fu Carlo Cecchini, che ritornava a Genova dopo aver completato gli studi in fisica alla Scuola Normale Superiore (SNS) di Pisa. Sempre in quel periodo Sandro stava dirigendo le tesi di laurea di due studentesse di Genova, Anna~Maria Mantero, purtroppo recentemente scomparsa, e Luisa Pedemonte, su temi di analisi di Fourier legati ad uno dei fondamenti degli insiemi lacunari, i polinomi di Rudin--Shapiro, e le funzioni di Rademacher e di Walsh viste come caratteri del gruppo di Cantor. In quel periodo si avvicinava all'analisi armonica anche Anna Zappa, che in seguito avrebbe scritto con Anna~Maria vari lavori molto profondi sulla trasformata di Poisson su gruppi liberi ed edifici di Bruhat--Tits.

Mi si lasci descrivere gli inizi della mia interazione con Sandro. All'epoca anche io studiavo Fisica alla SNS e cercavo di scrivere una tesi sulle rappresentazioni di alcuni gruppi di Lie, ma le applicazioni fisiche di questa teoria non mi entusiasmavano. Invece ero molto colpito dal fatto che Carlo, di cui ero concittadino ed amico, pur avendo una preparazione totalmente basata sulla fisica, riuscisse a sviluppare ricerche in matematica, delle quali mi parlava spesso con entusiasmo. Carlo mi suggerì di venire a Genova a parlare con Sandro, il quale mi propose di scrivere sotto la sua direzione una tesi su problemi di analisi armonica su gruppi compatti, un ambito le cui basi mi erano familiari perché ben conosciuto in fisica, ma i cui temi di ricerca erano certo molto diversi da quelli della fisica. Con un po' di titubanza, ma anche affascinato dall'entusiasmo di cui Sandro mi aveva fatto dono, accettai, completai nel corso dell'estate 1971 tutti i miei esami di Fisica e lasciai la SNS per scrivere la mia tesi di laurea a Genova. Ma prima ancora delle idee e l'entusiasmo che mi trasmise per la mia tesi, non posso omettere di menzionare un secondo dono che Sandro mi fece. Si occupava molto attivamente insieme a Carlo Pucci dell'insegnamento della matematica a livello avanzato, ma, poiché a quei tempi in Italia non era ancora stato istituito il dottorato di ricerca, come già accennato essi organizzavano corsi estivi di tipo dottorale tenuti da docenti stranieri, spesso statunitensi, i quali poi incoraggiavano i loro migliori studenti a frequentare scuole di dottorato nelle proprie università, scrivevano lettere di presentazione per aiutarli ad esservi ammessi e li aiutavano ad affrontare eventuali problemi iniziali di inserimento. Per questo fine crearono un ente autonomo, la SMI, di cui Sandro fu cofondatore ed in seguito Direttore e colonna portante, i cui corsi, a regime, si svolsero a Perugia. Ma nei primi due anni, nel 1970 e 1971, vennero tenuti invece a Pisa, nei locali della SNS, dove io ero ancora residente. Grazie a questo fatto e con l'aiuto e il consiglio di Sandro mi infiltrai abusivamente alle lezioni, e seguii ben tre corsi (gli iscritti regolari ne seguivano due). Pensavo che fosse un passo molto più lungo della gamba e che avrei fatto fiasco, ma non fu così, perché Sandro ed i bravi docenti mi incoraggiarono e mi diedero fiducia, ed io compresi a fondo gli argomenti avanzati dei tre corsi nonostante le mie conoscenze preliminari fossero solo quelle dei corsi di matematica per i fisici, esplorai nuovi universi mentali fatti di idee, ritrovai l'entusiasmo che avevo un po' perso studiando fisica: ma ora era entusiasmo per la matematica, che cambiò l'intera gamma dei miei interessi ed in nuce della mia futura creatività scientifica. Da quel corso estivo uscii entusiasta, ed affascinato dalla profondità di quei corsi di matematica. Sì certo, il mio esserne affascinato era naturale visto che all'inizio per me era tutto nuovo, ma la stessa opinione positiva la aveva avuta al corso estivo dell'anno prima Fulvio Ricci, di cui parleremo molto in seguito e che allora era uno studente di matematica alla SNS --- anche lui vi si era infiltrato prima di laurearsi, ed era uno studente bravissimo e molto competente. Conservo tuttora la mia opinione molto positiva dei corsi estivi SMI, che mi fa ammirare l'opera di Sandro nel crearli.

Fu da questo periodo che di fatto conobbi Sandro, quindi è da qui che questa presentazione ha una valenza personale affettiva --- ciò nondimeno continuerò per un po' a parlare anche di risultati scientifici. A quel corso estivo incontrai per la prima volta un giovane milanese, Paolo Maurizio Soardi, di cui in seguito parlerò molto, e lo incontrò per la prima vola anche Sandro: l'incontro tra i due fu, come vedremo, foriero di importanti sviluppi.

Subito dopo mi insediai a Genova, pur sempre iscritto a Fisica, nello stesso studio degli studenti locali di Sandro. Fu un'esperienza diversissima da quella che avevo avuto alla SNS. Il rapporto scientifico e personale di Sandro era costante e intenso con ciascuno di noi e sempre pieno di entusiasmo. Lui asseriva sempre di non sapere nulla di matematica, di non avere conoscenze ma solo idee, e quindi cercavamo di costruire ed esplorare insieme quell'universo di idee. Sandro aveva, e ha sempre avuto, una voce stentorea, con la quale ci chiamava nel suo studio, adiacente al nostro, attraverso la parete di cartongesso, quindi con scarso isolamento acustico, ma la cosa ci faceva sempre impressione, come una voce che calasse dal cielo. Allora andavamo nel suo studio e lì creavamo insieme idee --- spesso interrotti dalle continue telefonate di politica scientifica che Sandro riceveva, soprattutto da Carlo Pucci ed anche da Giunio Luzzatto ed altri collaboratori/amici. Sandro aveva un forte attaccamento con noi suoi allievi: passavamo a parlare con lui molto del tempo lasciato libero dalle lezioni universitarie e dalle telefonate di politica scientifica, ed al sabato mattina andavamo tutti insieme in piscina (dove non arrivavano le telefonate\dots). Poi, nel periodo di Natale, Irene ci invitava ai festeggiamenti a casa loro, dove preparava uno squisito \emph{eggnog}, che dopo mezzo secolo ho imparato io stesso a preparare, essendo diventato un cuoco migliore con il progredire dell'età, ed una volta l'ho portato per Natale a lei e a Sandro (una versione analcolica per loro, ma Irene assaporò con gusto anche qualche goccia della versione alcolica): era una affettuosa rievocazione del passato.

In tempi abbastanza brevi la mia tesi divenne un lavoro di ricerca su una classe di insiemi lacunari su gruppi compatti; non avrei mai pensato di essere in grado di dimostrare nuovi teoremi dopo solo pochi mesi dalla mia trasformazione in matematico, e senza il suo incoraggiamento probabilmente non sarebbe successo. In quel periodo egli ebbe come visitatori vari matematici, fra cui dall'Australia John Price, già studente di Robert~E. Edwards (maestro di Gaudry), insieme al quale scrisse due articoli \cites{19,20} sui polinomi di Rudin--Shapiro e su serie di Fourier aleatorie su gruppi compatti e moltiplicatori di $L^p$, e dalla Polonia Marek Bożejko, che molti anni dopo divenne un mio coautore. Questi temi erano nuovi nell'ambito di gruppi compatti non commutativi, ma in Italia non erano mai stati studiati neppure nell'ambito classico della retta reale o del cerchio unitario: Sandro stava quindi introducendo in Italia un nuovo ambito tematico, l'analisi armonica moderna, e simultaneamente stava costruendo dal nulla un gruppo di ricerca in questo campo, che rapidamente si espanse sempre più. In quel periodo vennero per una breve visita due giovani matematici milanesi, Paolo~Maurizio 
che abbiamo già citato e Leonede De~Michele, il quale, noto a tutti come Leo, sarebbe stato in seguito coautore di Sandro e suo collaboratore a lungo termine nelle iniziative di politica scientifica; il primo non scrisse mai articoli in collaborazione con Sandro --- ma Sandro apprezzava molto l'originalità, la competenza, la brillantezza dei suoi risultati. In seguito, entrambi ebbero una rilevanza fondamentale per l'impulso che diedero allo sviluppo del gruppo italiano di analisi armonica. Il rapporto con Sandro dei suoi allievi era quasi inevitabilmente e spesso indissolubilmente un rapporto di profonda amicizia: questa fu la leva che permise lo sviluppo iniziale di un gruppo di ricerca profondamente coeso. Quindi l'entusiasmo e la condivisione di intenti e di affetto che tutti noi avevamo nel far ricerca con Sandro furono, per me in particolare, assolutamente necessari. Ho già osservato come le tematiche e gli ambiti di ricerca di Sandro si stessero ramificando ed ampliando, dai gruppi abeliani degli inizi a gruppi compatti e poi a gruppi discreti non commutativi. Il \emph{trait d'union} che permetteva un trattamento unificato, più generale e astratto, si basava sulle C$^*$-algebre e sulle algebre di von~Neumann associate alle rappresentazioni dei gruppi localmente compatti, come presentate ad esempio nei testi ormai classici di Jacques Dixmier, a cui Sandro dedicò un'attenzione sempre crescente. Il suo articolo insieme al già citato Carlo 
\cite{22} riguardava appunto proiezioni di unicità per $L^p$ di un gruppo unimodulare non compatto definite in termini dell'algebra di von~Neumann generata dalle traslazioni su $L^2$: il nome ``unicità'' incarnava un loro risultato sorprendente, ossia che per ogni $1\leq p<2$ esiste un sottospazio chiuso non banale di $L^2$ invariante per traslazioni destre che non contiene funzioni di $L^p$ oltre alla funzione nulla.

Farò qualche cenno ora, e più volte in seguito, al mio lavoro di ricerca con Sandro. In quel periodo, estate 1972, avevo completato il mio risultato sugli insiemi lacunari su gruppi compatti, ed allora mi indirizzò allo studio di un ambito completamente opposto, il gruppo libero, ossia un gruppo discreto infinito in qualche senso ``universale'', che molti anni dopo divenne per me l'ambito cardine della mia ricerca. Sandro aveva un'abilità straordinaria per passare da un ambiente all'altro e riformulare conseguentemente problemi e metodi, ma io no, e ne ero un po' confuso. Aveva in mente due risultati recenti legati ai gruppi liberi. Il primo era dovuto a Michael Leinert, un comune amico purtroppo recentemente scomparso, il quale aveva scoperto che un sottoinsieme libero $E$ di un gruppo discreto ha la proprietà che, per un'opportuna costante $C$, ogni funzione definita su $E$ ha norma di von~Neumann maggiorata da $C$ per la norma $L^2$. Il secondo risultato, dimostrato qualche anno dopo, era una fondamentale disuguaglianza dovuta a Uffe Haagerup fra la norma $L^2(\setF)$ e la norma $A(\setF)$ per le funzioni con supporto sulle parole di lunghezza $n$ su un gruppo libero $\setF$, dal quale Sandro traeva importanti conclusioni sulle funzioni in $B(\setF)$ (lo spazio delle combinazioni lineari delle funzioni di tipo positivo) e sui moltiplicatori di $A(\setF)$ (e qualche anno dopo ne trasse la caratterizzazione del preduale dello spazio dei moltiplicatori di $A(\setF)$ \cite{30}). Questo secondo argomento per noi fu importante per l'ultimo capitolo di un libro che lui ed io scrivemmo dieci anni dopo \cite{34}, nel quale ormai consideravamo il gruppo libero da un punto di vista geometrico, come un gruppo semplicemente transitivo di automorfismi di alberi omogenei. Il primo invece si tradusse nel nostro primo articolo congiunto \cite{21}, dove dimostrammo per $\setF$ un altro fatto che in analisi armonica classica è chiaramente falso, ossia che esistono moltiplicatori di $A(\setF)$ che non appartengono a $B(\setF)$. In quell'estate, oltre che presentare questo risultato nella mia prima conferenza scientifica internazionale, frequentai questa volta in veste ufficiale i corsi estivi matematici e, grazie alle lettere di presentazione dei miei docenti e di Sandro, fui ammesso come studente di dottorato alla University of Maryland, dove simultaneamente egli era stato invitato come professore visitatore per un anno. Subito prima era nato il suo secondo figlio, Niccolò. Al termine di quell'anno Sandro ritornò a Genova e suggerì problemi scientifici a Giancarlo Mauceri --- che tutti abbiamo sempre chiamato Gianni --- ed in seguito ne divenne coautore in~\ocite{28}, articolo che utilizza risultati ottenuti da Sandro in~\ocite{25}. Gianni ed io ci conoscevamo un pochino sin dall'inizio della mia permanenza a Genova, ma poi nacque una profonda amicizia quando, durante un mio breve ritorno a Genova dall'America, studiammo insieme per prepararci ad un concorso per diventare Assistenti Universitari: la collaborazione in matematica, e presumo in ogni altro campo, è una condivisione di idee ed aspettative che spesso porta a reciproco apprezzamento ed amicizia. Dopo il mio ritorno definitivo dall'America, Gianni ed io stabilimmo una proficua collaborazione, da cui scaturì la costruzione di esempi interessanti di gruppi non compatti la cui algebra di von~Neumann è puramente atomica: questi argomenti erano basati su strumenti di teoria delle rappresentazioni che avevo appreso nei miei studi a Maryland e su altre tematiche legate a quelle che egli aveva studiato insieme a Sandro. Naturalmente continuò e si intensificò la nostra amicizia, che dura tuttora: ricordo quando, con le rispettive consorti, anche loro amicissime, avevamo trovato un posto in riva ad un ruscello dove fare campeggio libero (nostre le due uniche tende) e lì, sedendo sull'erba, Gianni ed io oltre a chiacchierare e scherzare facevamo anche ricerca insieme.

In quegli anni Sandro pubblicò vari risultati interessanti su moltiplicatori su $L^p(G)$ (con $G$ gruppo compatto) che tendono a zero all'infinito sul duale di $G$ \cite{23}, e sull'esistenza di funzioni in $B(G)$ che tendono a zero all'infinito ma non appartengono ad $A(G)$ per una vasta classe di gruppi unimodulari non compatti, quelli la cui algebra di von~Neumann non è puramente atomica \cite{24}. Al mio ritorno definitivo dall'America io trasferii la mia posizione di assistente universitario a Perugia, seguendo Sandro che nel frattempo, nel 1974, vi si era trasferito come professore ordinario. A causa della mia permanenza in America, non avevo avuto contatti diretti con Sandro per un paio d'anni e non posso parlare di quel periodo per testimonianza diretta, ma è importante menzionare che proprio allora venne a Genova un altro giovane matematico australiano originariamente studente di Gaudry, Michael Cowling, che in seguito risultò fondamentale per lo sviluppo dell'analisi armonica in Italia (e non solo). In quegli stessi anni venne coinvolto nel gruppo di analisi armonica anche Fulvio Ricci --- ho già riferito di averlo conosciuto a Pisa quando eravamo entrambi studenti alla SNS: egli era stato studente di Edoardo Vesentini --- eminente studioso di mappe olomorfe, spazi analitici e varietà di Kähler, funzioni analitiche in più variabili complesse ed analisi funzionale, matematico di primaria importanza. Edoardo Vesentini fu grande amico di Sandro e di Carlo Pucci e compartecipe degli stessi aneliti verso un sistema più internazionale di studi a livello del dottorato di ricerca, ed anche per perseguire questi obiettivi fu in seguito Senatore della Repubblica e Vicepresidente della Commissione per la pubblica istruzione ed i beni culturali. Fulvio aveva conosciuto Sandro a Pisa nell'estate del 1970, alla prima edizione dei corsi estivi matematici della SMI, e poi, indirizzato da Vesentini, aveva studiato analisi di Fourier nel suo dottorato all'Università del Maryland sotto la guida di John~J. Benedetto. In seguito Fulvio ebbe un ruolo apicale nello sviluppo di questo gruppo di ricerca a Pisa, poi a Torino, poi nuovamente a Pisa, così come Gianni a Genova, Paolo~Maurizio (e qualche anno dopo Leonardo Colzani e Giancarlo Travaglini) a Milano, Michael prima a Genova e poi in Australia, ed in seguito Wolfgang Woess prima a Milano e poi a Graz; molti dei loro allievi continuarono questo sviluppo avviando alla ricerca numerosi altri giovani valenti matematici. Questo processo prosegue ancora, generazione dopo generazione; mi scuso del fatto che me ne sia impossibile un'analisi individuale. Per lo sviluppo (e la coesione) del gruppo è stata di valore inestimabile anche l'attività infaticabile di Leo.

A Perugia Sandro coprì una cattedra di Matematica generale alla facoltà di Economia, e lì nacque il suo terzo figlio Lorenzo. In quella sede Sandro diresse gli inizi della ricerca di Mauro Pagliacci, il quale in seguito diede numerosi contributi all'analisi armonica su gruppi che agiscono su alberi omogenei e poi ne trasse applicazioni alla matematica finanziaria e divenne Professore in questa disciplina. A Perugia venivano periodicamente a parlare di matematica con Sandro i sunnominati colleghi milanesi Paolo Maurizio e Leo, e sempre da Milano il compianto Alberto Alesina. Nel 1976, Sandro contribuì significativamente alla diffusione ed ai legami internazionali dell'analisi armonica in Italia organizzando nella sede dell'INdAM, a Roma, un congresso di livello eccelso in Analisi Armonica (fu il primo congresso da lui organizzato, ne seguirono poi altri a Cortona, di uguale successo: gli piaceva molto organizzarli). A quel congresso parteciparono molti matematici stranieri che in seguito furono determinanti per lo sviluppo della ricerca italiana in questo campo: ad esempio Guido~L. Weiss, ma anche Michael che in quel periodo stava completando un fondamentale risultato sul fenomeno di Kunze e Stein --- e qui voglio rammentare un altro conferenziere, Pierre Eymard, perché, come vedremo, la sua presentazione fu determinante per la ricerca di Sandro e mia di molti anni dopo. Nell'estate 1976 era stato trasferito a Perugia uno dei corsi estivi matematici della SMI, che di solito si svolgevano a Cortona: si trattava proprio di un corso di analisi armonica, tenuto da Ray~A. Kunze e Guido Weiss. A partire da quel momento Guido ebbe un ruolo determinante nel fornire appoggio e idee per lo sviluppo dei giovani analisti armonici italiani, molti dei quali infatti studiarono alla Washington University a Saint~Louis, ed egli stesso venne più volte in visita in Italia: fu, dopo Sandro, uno dei fulcri intorno a cui si evolse questo gruppo di analisi armonica. Furono molto influenzati da Guido, in particolare e fra tanti altri, Gianni e Fulvio, il quale in seguito ne divenne anche autorevole coautore. In quegli anni il gruppo di ricerca in analisi armonica prendeva corpo con un notevole affiatamento. Tenevamo settimanalmente o bimensilmente un seminario, che Sandro ci aveva esortato a svolgere allo scopo di imparare da noi argomenti che avevamo studiato, alcuni di noi all'estero. Al seminario partecipavano i vari colleghi ed amici di Genova (Michael ed Anna~Maria), di Pisa (Fulvio), Milano (Paolo~Maurizio e Leo), Perugia (Walter Betori e Mauro) ed in seguito, dopo che vi si trasferì Sandro, anche Roma; avevamo scelto come punto d'incontro Firenze, dove Sandro comunque aveva piacere di incontrare Carlo Pucci ai vicini uffici del CNR al fine di concertare azioni nell'interesse della matematica italiana e dei suoi giovani ricercatori. Anzi meglio: Pucci sapeva che Sandro stava per un paio d'ore al seminario lì a due passi, ed approfittava di questa ghiotta opportunità per chiamarlo al telefono ogni quarto d'ora circa per concordare iniziative politiche, chiedendogli di raggiungerlo subito\dots\ ed ahinoi, il forte senso del pubblico interesse di Sandro spesso era altrettanto forte della sua passione per la ricerca scientifica, ed alla fine il suo senso del dovere prevaleva e, rattristato ed insofferente del destino avverso, ci lasciava per un po' --- e noi continuavamo da soli, altrettanto rattristati perché volevamo condividere con lui le idee matematiche, che al momento trovavamo più importanti di quelle organizzative, sebbene sapessimo che stava organizzando iniziative a favore dei giovani matematici come noi. Nel tenere questi seminari, miravamo anche a condividere le nostre conoscenze in vari settori dell'analisi armonica, perché i periodi trascorsi all'estero da molti di noi ci avevano fornito competenze diverse. In effetti, ciascuno dedicava un grande impegno ad ampliare le proprie. Ad esempio, Gianni aveva passato il biennio 1976--1978 negli Stati Uniti, il primo dei quali insieme ad Alberto Alesina e Leo (ed anche in tal caso nacque una profonda amicizia), ed aveva cominciato a interessarsi agli operatori di \emph{twisted convolution}. Al suo ritorno dall'America, Gianni, Fulvio ed io ci incontravamo periodicamente a Pisa per collaborare affrontando questa tematica per me assolutamente nuova: erano per me momenti fondamentali, apprendevo da due amici bravissimi e già esperti e scoprivo con gioia nuovi orizzonti. Analoghi contatti su altri temi avevano frequentemente luogo fra gli analisti armonici attivi a Milano e a Genova, ed in seguito a Roma. Eravamo un gruppo di amici in continuo movimento!

A proposito, il parlare di Firenze mi rammenta un altro carissimo amico di Sandro (e poi anche mio), che era stato allievo di Carlo Pucci ed era professore lì: Giorgio Talenti. Giorgio si occupava di equazioni differenziali alle derivate parziali, ma con una cultura molto vasta che quindi spaziava largamente anche nell'analisi di Fourier, e con un gusto sopraffino. Oltre a riconoscere la bellezza nella matematica, la riconosceva negli animi. Si interessò molto, insieme a Stefano Campi, alla visualizzazione della matematica a scopi didattici, e nel periodo in cui anche io me ne occupavo e scrivevo programmi di grafica per computer, tutti e tre insieme scrivemmo un libro su esperimenti numerici e grafici al computer per gli studenti universitari che dovevano apprendere il calcolo in una o più variabili: erano i tempi in cui ancora non esistevano applicativi specifici di questa natura. A quell'epoca ero alla Sapienza, e Sandro mi guardava con grande rispetto mentre creavo grafici tridimensionali al computer: credo che equiparasse il computer alla magia.

Non sorprende che Sandro, vista la sua personalità radiosa, avesse tanti amici: io sono consapevole quasi soltanto di quelli che erano matematici. Normalmente, Sandro incontrava i matematici ai congressi o durante le visite scientifiche. A volte, a questo tramite se ne aggiungevano altri. Mi si lasci citare il caso di Nicholas Varopoulos, un famoso matematico greco residente a Parigi, a cui mi riferirò come Nick: trascorreva una parte considerevole delle sue vacanze in una casa nell'isola greca di Corfù. Anche Sandro e Irene trascorrevano abitualmente le vacanze nella loro casa a Corfù, che Irene aveva ereditato da sua madre. Ebbene, l'amicizia fra Sandro e Nick fu consolidata dalle frequentazioni estive fra le onde del mar Ionio --- oltre che dalle loro ricerche scientifiche di rilievo.

Sandro era molto legato alla Grecia, la patria di sua moglie Irene. Avendo studiato la storia e la letteratura della antichità greca al liceo, aveva visitato con i figli e nipoti tutti i luoghi archeologici della Grecia, e visto rappresentazioni delle tragedie greche in teatri antichi. Era anche cultore della storia bizantina e visitò due volte Monte Athos, portandosi dietro, lungo l'interminabile sequenza di gradini, figli e nipoti maschi, dato che le femmine non sono ammesse. Era orgoglioso di queste due visite e le raccontava spesso, persino negli ultimi anni di vita.

Anche il loro primogenito Giovanni trascorreva le vacanze in una isola ionia vicino a Corfù, e veleggiava per l'arcipelago con la sua barca a vela. Sandro era molto compiaciuto delle abilità velistiche di Giovanni, come di tante altre: diceva che Giovanni sapeva fare tutto, e che in questo assomigliava a me --- in qualche modo ciò mi faceva sentire assimilato a figlio, come asseriva Irene.

La Grecia, ed in particolare Corfù, ebbero molta importanza nella sfera emotiva di Sandro e Irene, quindi mi si lasci narrare qualche memoria che ho di loro in quell'isola. Sentire Sandro e Irene parlare di Corfù, di persone celebri che vi erano originarie, dell'uliveto Rothschild, del fatto che metà dei maschi locali si chiamasse Spiridione come il santo patrono, era davvero divertente.
Nell'estate 2007 la mia famiglia ed io andammo a trascorrere le vacanze nella loro villa. Adoravano quell'isola, e ricordo come ci trasmisero il loro entusiasmo --- ricordo anche la battaglia che combattemmo contro il formicaio delle termiti che si stavano rosicchiando la trave di legno che sosteneva il piano superiore: la casa restò in piedi, quindi penso che vincemmo quella battaglia, anche se forse non la guerra.

L'anno successivo Sandro avrebbe celebrato il suo settantesimo compleanno, ed io avevo intenzione di organizzare a Corfù il congresso annuale di analisi armonica (di cui parlerò ampiamente in seguito). La sede più appropriata sarebbe stata il Rettorato dell'Università dello Ionio, Palazzo Grimani, prestigioso edificio che era appartenuto all’Accademia Ionica fondata a Corfù nel 1824 e dove fin da allora l'Accademia aveva organizzato i primi studi universitari ellenici, fino a quando fu chiusa all'inizio della sovranità greca nel 1865. E così Sandro ed io, accompagnati dalle rispettive famiglie, andammo in visita cerimoniale dal Rettore e dalla Contessa Bulgari di Corfù, la quale aveva un ruolo preminente nelle questioni culturali ed organizzative, in particolare quelle nella tradizione della antica Accademia Ionica.

Per spiegare meglio il contesto di quella visita, occorre delineare due punti: l'aristocrazia dell'isola ed il modo scherzoso in cui la descrivevano Sandro ed Irene, il cui senso dello spirito affiorava sempre e ovunque.

I Bulgari di Corfù sono una famiglia dell’alta aristocrazia che si era stabilita fin dal quindicesimo secolo nell'isola, che allora era possedimento della Repubblica di Venezia. Corfù non fu mai ottomana, rimase possesso veneziano fino alla caduta della Serenissima sotto Napoleone, dopo il 1815 fu protettorato britannico. La famiglia Bulgari di Corfù è stata per lungo tempo custode delle reliquie di San Spiridione, santo patrono molto venerato nella Chiesa ortodossa di Corfù, dove viene onorato ben cinque volte all’anno. Si deve a questa famiglia la costruzione, secoli fa, della chiesa di San Spiridione, e quindi la tutela spirituale della popolazione dell'isola, che mantennero coesa nella difesa dagli assedi ottomani.

Con la loro consueta, accattivante miscela di spirito e arguzia, Sandro e Irene mi spiegarono che la Contessa ci teneva a precisare di non avere alcuna relazione con la famosa famiglia Bulgari di gioiellieri di Roma, che definiva ``i commercianti di Roma'' (venditori ambulanti al Pincio\dots\ indicazione peraltro sostenuta dal documento storico della ditta Bvlgari ``The Life of Sotirio Bvlgari, the Founder of Bvlgari'' su YouTube) --- si tratta dei discendenti di quel membro di un ramo distante della famiglia, Sotirios Voulgaris (Sotirio Bulgari), nato nel 1857 vicino a Ioannina, nell’Epiro allora ottomano, dotato di indiscutibile maestria artigianale nel produrre gioielli d'argento (la lavorazione della filigrana d’argento era ispirata sia dalla tradizione ottomana sia da quella greca). Sotirio Bulgari si era trasferito a Napoli nel 1881 e l'anno dopo a Roma, dove aveva cominciato a vendere i suoi monili d'argento al Pincio e subito dopo aveva lavorato per importanti gioiellerie ed infine aveva aperto negozi di successo prima a Roma, poi in Italia ed all'estero, e fondato un marchio che divenne in pochi decenni noto in tutto il mondo come sinonimo di lusso: il disdegno della contessa sottolineava il fatto che Sotirio Bulgari non era detentore del titolo nobiliare!

Possiamo ora tornare al racconto di quella visita al palazzo Grimani, diventato ormai quasi un museo che rievoca le vecchie glorie dell'Accademia Ionica, con ritratti alle pareti ed una ``aula magna'' che, in contrasto con il nome, è una minuscola saletta della quale intendevo chiedere l'uso per organizzare il Congresso celebrativo del settantesimo compleanno di Sandro. Lui ed io fummo guidati ad esplorare ogni stanza del palazzo e dedicammo ossequiosa attenzione a tutti i ritratti dei professori di quasi due secoli prima, ed io sorridevo fra me e me all'idea di celebrare Sandro in mezzo ai ritratti di tutti quei vecchi baroni (beh, forse alcuni erano conti\dots). Poiché Sandro trascorreva a Corfù ogni estate e di fatto era considerato uno del posto, parlava perfettamente greco ed aveva moglie greca, e fu facile spiegare la rilevanza scientifica di un evento culturale ``storico'' quale il Congresso di Analisi Armonica, che ormai stava raggiungendo la trentesima edizione, e l'appropriatezza di svolgerlo in quel prestigioso edificio. Presentando Sandro come famoso corfiota, in tutto degno degli onori della (ormai estinta) Accademia Ionica, ed evitando accuratamente nel colloquio con la Contessa di menzionare Roma ed i suoi commercianti, ebbi pieno successo nella mia missione diplomatica: l'auletta ci fu generosamente concessa, ad un costo contenuto. Ma poi emerse che i fondi che avevo richiesto ed ottenuto per il Congresso non potevano essere usati per finanziare visite di scienziati stranieri in località non italiane. Perciò il destino fu avverso al sogno dell'edizione Ionica del nostro congresso.

Meglio riassumere fin da subito l'evoluzione del gruppo di analisi armonica dopo il biennio che Sandro passò a Perugia; mi scuso per la lunga divagazione che spezza la sequenza cronologica. Si aggiungevano al gruppo altri giovani: ad esempio, a Milano Giancarlo e Leonardo già menzionati prima, nonché Saverio Giulini, di cui compiangiamo la scomparsa, che in seguito si trasferì a Genova; in una successiva ``generazione'' Maria~Gabriella Kuhn ed Elena Prestini (il cui ingresso in analisi armonica era stato indipendente da Sandro, avendo scritto la tesi di dottorato a Princeton sotto la guida di Charles~L. Fefferman), la quale poi fu a Milano e molti anni dopo venne a Roma ``Tor Vergata''; nello stesso periodo da Trento a Roma venne anche Welleda Baldoni, che con Sandro aveva scritto la tesi di laurea a Genova e dopo il suo dottorato americano sotto la guida di Nolan~R. Wallach si dedicò alla teoria delle rappresentazioni dei gruppi di Lie semisemplici. E via via tanti altri, così numerosi che è impossibile elencarne tutti i nomi,
ma mi si lasci almeno tracciare un succinto quadro storico, sebbene gravemente incompleto, della evoluzione successiva. Vorrei rammentare Stefano Meda, che scrisse la tesi con Gianni, e dopo scrisse anche qualche articolo su spazi di Lorentz su alberi, quindi nell'ambito iniziato da Sandro e continuato da Michael, ma poi tanti altri articoli sull'operatore di Ornstein--Uhlenbeck ed il nucleo del calore, anche in collaborazione con Alberto Setti (entrambi erano a Milano, poi Alberto si trasferì a Como), e con Michael e soprattutto con Gianni. A poco a poco, il gruppo si allargò sempre di più, ormai con interessi molto variegati e via via meno vicini a quelli di Sandro; qui rammento ad esempio, sempre a Milano Luigi Fontana e Carlo Morpurgo, anche loro poi studenti a Washington University (dopo il dottorato ritornarono a Milano e poi Carlo andò alla Università del Missouri), ed Enrico Laeng che fu studente a Roma, poi a Yale e dopo il dottorato andò a Milano, e Filippo De~Mari, che scrisse la tesi di laurea con Gianni a Genova e fu studente di dottorato a Washington University, e dopo ritornò a Genova. Con piacere rammento qui uno dei miei coautori più frequenti, Enrico Casadio~Tarabusi, che, dopo aver completato gli studi dottorali alla SNS sotto la guida di Vesentini, fu per quasi un anno post-doc a Maryland, dove insieme a Joel ed a Carlos Berenstein scrivemmo il nostro primo articolo su trasformate di Radon su alberi completamente generali, e poi fu a Trento e dopo a Roma ``La Sapienza'' (egli fu l'ultimo coautore di Sandro e verrà citato molte volte in seguito come Enrico nonostante la omonimia con Enrico Laeng). Ed ancora, a Milano Luca Brandolini e Maura Salvatori, a Genova Ernesto De Vito, a Torino Anita Tabacco, che, dopo essersi interessata a lungo di \emph{wavelets}, (come, sempre a Torino, fece Elena Cordero) in questi ultimi anni ha rivolto l'attenzione a problemi di analisi su alberi, cosa che avrebbe interessato Sandro se in tali anni si fosse occupato ancora di matematica, e vari anni dopo sempre a Torino anche Marco Peloso, con una formazione più basata sull'analisi complessa. Ancora a Milano, rammento Marco Vignati ed anni dopo Giacomo Gigante, ed a Bologna Nicola Arcozzi, di cui mi onoro di essere stato coautore insieme ad Enrico ed a Fausto Di Biase (che nel seguito di queste pagine diventerà una presenza frequente) anche se non proprio in analisi armonica bensì in teoria del potenziale, campo in cui egli era già un ricercatore molto autorevole: in questo contesto Fausto ci guidò a sviluppare con successo una sua idea molto acuta ed originale su un problema classico. Nicola aveva scritto la tesi di laurea a Milano sotto la guida di Leo, e poi era stato studente di dottorato di Albert Baernstein II a Saint~Louis; ha avuto un ruolo importante nell'espansione del gruppo di analisi armonica, guidando numerosi allievi su problemi fini non solo di analisi armonca (disuguaglianze di Hardy, trasformate di Riesz, spazi di Besov, spazio di Dirichlet, teoremi di Paley--Wiener, misure di Carleson), ma anche di teoria del potenziale in vari ambiti, inclusi gli alberi. E cito anche Sandra Saliani (della quale scriverò più avanti) a Potenza e più tardi a Napoli, anche lei almeno inizialmente studiosa delle \emph{wavelets} (Sandro ed altri le chiamavano ``ondine'', io ``ondicelle'', ed ogni tanto discutevamo sulla appropriatezza di queste scelte. Io asserivo che le ondine fossero creature della mitologia tedesca introdotte con questo nome da Paracelso e noi dovessimo rispettarne il diritto di autore. Per sciogliere il dilemma qui uso l'espressione canonica anche se è inglese). E poi ancora, anni dopo, Maria Vallarino a Genova, che scrisse la tesi sotto la guida di Stefano e poi si trasferì a Torino e studiò trasformate di Riesz, spazi di Bergman, spazi di Damek--Ricci e spazi di Hardy, e che negli ultimi anni ha anche lei allargato il suo orizzonte all'analisi armonica su alberi. Ed anche Vittoria Pierfelice, proveniente dall'ambito delle equazioni alle derivate parziali e che poi studiò con Fulvio a Pisa ed in seguito si trasferì a Orleans. E Bianca Di~Blasio che da Milano venne a Roma ``Tor~Vergata'' e poi andò a Milano Bicocca, e Francesca Astengo a Genova: lavorarono insieme, sotto l'ispirazione di Fulvio nonostante le sedi fossero differenti, e dopo anche in collaborazione con Michael, su problemi sulla trasformata di Fourier nel senso di Helgason su vari gruppi di Lie,
trasformate di Cayley su gruppi $AN$ e coppie di Gelfand su gruppi di Heisenberg;
Francesca scrisse anche un articolo di analisi su alberi. Ci furono anche Paolo Ciatti e Valentina Casarino, a Padova, che ottennero risultati su aspetti geometrici dell'analisi armonica, in particolare per le proprietà di limitatezza di certi operatori di proiezione (proiezioni spettrali del Laplaciano su gruppi di Lie nilpotenti e proiezioni armoniche su sfere complesse e quaternioniche), e moltiplicatori spettrali per un operatore di Ornstein--Uhlenbeck non simmetrico. Ed a Milano ci fu Mauro Maggioni, un altro \emph{alumnus} di Washington University, il compianto Roberto Camporesi a Torino, e moltissimi altri, una lista troppo lunga da elencare e che probabilmente si va espandendo anche ora mentre scrivo: delle ultime generazioni mi vengono in mente, perché ho ascoltato le loro conferenze ai congressi annuali di analisi armonica di cui sto per parlare, i nomi di Andrea Carbonaro, Fabio Nicola, Alessio Martini, Matteo Monti, Tommaso Bruno, Alessandro Ottazzi, Federico Santagati, Matteo Levi, Bianca Gariboldi, Alessandro Monguzzi, Mattia Calzi, Stefano Vigogna, 
Francesca Bartolucci, ed altri almeno inizialmente studiosi di analisi armonica, come Rita Pini ed Andrea Calogero, ma chissà quanti altri; qui mi fermo, anche perché Sandro non ha neppure conosciuto parecchi di questi giovani delle ultime generazioni.

In effetti, a Washington University, nell'atmosfera accogliente e stimolante creata da Guido ed anche da Ronald Coifman, Richard Rochberg, Albert Baernstein II, Gary Jensen, Steven Krantz, Edward Wilson, Mitchell Taibleson, ciascuno dei quali fu direttore di tesi o collaboratore di ricerca di matematici italiani, si era creata una tradizione di studenti e docenti italiani (ed anche spagnoli), non tutti in analisi armonica, ma tutti incoraggiati almeno in parte dalla spinta di Sandro alla internazionalizzazione degli studi ed alla creazione di iniziative di supporto finanziario statale. Il modo affettuoso e compartecipe, quasi paterno, con cui gli studenti del Dipartimento di Matematica di Washington University erano accolti ed assistiti era simile al rapporto che Sandro stesso aveva con i propri. Molti dei matematici che ho citato nelle pagine precedenti studiarono lì oppure vi trascorsero dei periodi di ricerca o insegnamento, a partire da Fulvio e Gianni, ed in seguito tanti altri (inclusi Sandro e me). Vari giovani italiani studiarono a Saint~Louis, molti analisti armonici che abbiamo nominati prima, ma anche altri che maturarono interessi diversi, ad esempio di geometria come Marco Rigoli, che poi partecipò attivamente a vari congressi annuali italiani di Analisi Armonica. E naturalmente qui non sono ancora nominati tanti altri che si formarono o completarono la loro formazione a Roma: ad essi accennerò in seguito. Inoltre, non posso dimenticare il ruolo fondamentale che un altro amico di Sandro e di tutti noi, Peter Sjögren, di Göteborg ma con molti legami con l'Italia, svolse per la formazione di tanti giovani e la collaborazione con tanti altri, giovani o meno, incluso me. Nei quaranta anni in cui Sandro sviluppò in Italia l'analisi armonica vi furono alcuni altri matematici italiani attivi in questo campo ma non direttamente membri di questo gruppo e spesso non residenti in Italia, ed altri che furono attivi in altri campi ma almeno marginalmente anche in analisi armonica; vari di loro parteciparono occasionalmente o periodicamente ai congressi annuali del gruppo. Mi dispiace non nominarli esplicitamente, ma non posso essere certo di ricordarli tutti, ed ometterne solo alcuni sarebbe scortese.

Ho già implicitamente accennato al fatto che sto per esporre, ossia che l'ambito geometrico delle ricerche di Sandro stesse a poco a poco per concentrarsi su alberi omogenei. Anticipo quello che spiegherò poi in dettaglio, ossia che sviluppò importanti teorie sulle rappresentazioni di gruppi che agiscono su questi alberi e su fenomeni di diffusione e passeggiate casuali su di essi. Il trasferire ad alberi teorie che sono consolidate in contesti più classici, come gli spazi euclidei ed i gruppi di Lie, può essere considerato un settore di nicchia nella ricerca scientifica, ancorché abbia legami amplissimi con altri settori cruciali della matematica, come la teoria dei grafi (non solo alberi omogenei), la teoria discreta del potenziale, il comportamento asintotico delle autofunzioni dell'operatore di Laplace discreto e di altri operatori di transizione, e quindi le passeggiate casuali e le catene di Markov, tutti campi nei quali vari collaboratori di Sandro hanno svolto indagini di successo. Inevitabilmente, e non solo in Italia, dall'ambiente accademico consolidato le ricerche di nicchia in ambiti esotici sono inizialmente considerate con un certo snobismo (se mi si permette l'anglicismo), e questo un po' è accaduto anche fra i membri del gruppo di analisi armonica, pur senza ridurre la stima reciproca. Il fatto che ho osservato più sopra riguardo alle linee di ricerca di alcuni degli analisti armonici delle generazioni successive (ma che è vero anche per vari altri le cui linee di ricerca non ho avuto spazio per descrivere), ossia che dopo molti anni abbiano cominciato a rivolgersi anche allo studio dell'analisi su alberi, sottolinea in modo significativo il ruolo apicale e preveggente che Sandro ha avuto nello sviluppo scientifico del gruppo.

Ritorniamo al decennio in cui conobbi Sandro. Proprio visto il fatto che all'interno del gruppo la collaborazione scientifica stava diventando un legame fra amici in sedi lontane, e si poteva perdere traccia di quali risultati gli altri stessero cercando, nel 1978 Sandro varò un'iniziativa scientifica illuminata: un convegno-seminario annuale informale, in cui ci incontrassimo tutti non per farci vanto dei nostri risultati presentandoli in modo succinto (quindi spesso in maniera poco comprensibile), bensì per far capire a tutti i dettagli dei problemi che ostacolavano il nostro progresso e le loro inerenti difficoltà, con l'auspicio che altri di noi fornissero spunti per risolverle. Quel primo convegno fu davvero informale, quasi da festa paesana, come un incontro al bar di un gruppetto di amici. Michael ed Anna Maria trovarono una pensione sul costone d'una montagna a Sant'Ilario, vicino a Genova, e presero in affitto per una settimana una stanzetta in cui portarono dall'università una lavagna su ruote con la lastra di ardesia rotante per poter scrivere anche sul lato posteriore. Ma per farlo bisognava prima scostare la lavagna dal muro, perché le piccole dimensioni della stanza lasciavano solo un ristretto passaggio dietro, nel quale doveva comunque strisciare, durante le conferenze, chi aveva necessità di andare al bagno, che era proprio lì dietro\dots\ Tutti noi siamo adesso orgogliosi che dopo quasi mezzo secolo quel seminario con cadenza annuale si svolga ancora, sebbene forse ed inevitabilmente non più con lo stesso spirito pionieristico.

In quel periodo, fra il 1976 ed il 1978, io stavo anche completando la mia tesi dottorale, dove dimostrai che un gruppo di Lie non compatto non può avere algebra di von~Neumann puramente atomica. Ma allo stesso tempo, grazie alle conoscenze sulle rappresentazioni del gruppo di Heisenberg maturate nei miei studi americani ed alle profonde conoscenze di Sandro sulle funzioni che operano sull'algebra $B(G)$, nell'anno passato insieme a Perugia avevamo dimostrato che, se $G$ è un gruppo con centro non compatto oppure un gruppo nilpotente, solo le funzioni intere operano sulla sottoalgebra $B_0(G)$ delle funzioni in $B(G)$ che tendono a zero all'infinito. Subito dopo questo lavoro e quel convegno informale, a fine 1976, Sandro si trasferì all'Università di Roma, ritornando come Professore Ordinario nell'istituto in cui aveva cominciato i propri studi universitari. Un anno dopo, con molte preoccupazioni alla prospettiva di dover trovare casa in una città così difficile e complessa ed a me del tutto ignota, mi ci trasferii anch'io. Restai all'Università di Roma per 10 anni, dal 1977 a fine 1987, e quella fu l'ultima delle quattro sedi universitarie che ho condiviso con Sandro; lì scrivemmo i nostri ultimi lavori congiunti, che furono alcuni dei più significativi della mia vita scientifica e i più importanti della nostra collaborazione. Prima che questo accadesse, però, conobbi un altro suo studente, il primo che Sandro ebbe a Roma, Claudio Nebbia, un matematico di grande acume con il quale in seguito Sandro scrisse un bellissimo libro \cite{39} con risultati profondi in analisi armonica su alberi. Devo anche assolutamente citare Joel Cohen, che, dopo aver incontrato Sandro quand'era studente al Massachusetts Institute of Technology (MIT), dove avevano uffici adiacenti, in quel periodo venne quasi annualmente in Italia e visitò ogni volta Roma. Joel divenne non solo analista armonico ma anche grande amico e coautore di Sandro \cite{38}, cittadino italiano e professore a Bari, dove avviò all'analisi armonica Sandra Saliani (che incoraggiò ad andare all'Università del Maryland come studentessa di dottorato di John J. Benedetto) ed in parte Silvia Romanelli, con la quale organizzò nel 1985 l'edizione pugliese del convegno annuale. La prima visita di Joel in Italia fu ai corsi estivi SMI di Perugia nel 1977, dove conobbe Flavia Colonna, che in seguito sposò. Flavia lo seguì all'Università del Maryland, dove fu studentessa di dottorato; scrisse la propria tesi in analisi complessa sotto la guida di Maurice H.~Heins, e poi ebbe una produzione copiosa di eccellenti articoli su vari temi, fra cui, parecchio tempo dopo la tesi, anche vari articoli in analisi armonica, funzioni armoniche e teoria del potenziale su alberi congiuntamente con Joel, Enrico, me ed altri coautori.

Insieme a Claudio ed a me, Sandro organizzava seminari periodici rivolti anche agli studenti di tesi. In quei tempi, alla fine del decennio dal 1970 al 1980, i suoi temi di ricerca preferiti non erano più soprattutto inerenti alla teoria dei moltiplicatori di $L^p$ su vari gruppi localmente compatti ed alle funzioni di tipo positivo, anche se scrisse un bell'articolo con Leo \cite{32} dove costruirono su un gruppo libero $\setF$ funzioni in $B(\setF)$ con valori sui generatori preassegnati arbitrariamente (naturalmente di modulo non superiore ad uno). Leo veniva spesso a Roma per completare quella ricerca con Sandro; lo rammento con piacere perché già allora eravamo grandi amici. Le linee di ricerca di Sandro si rivolgevano sempre di più all'algebra di von~Neumann ed alla C$^*$-algebra dei gruppi, ma con particolare interesse per il gruppo libero. Studiare insieme nel corso di un seminario può produrre utili conseguenze: Sandro ed io presentammo un articolo espositivo di Hillel Furstenberg, \textit{Random walks and discrete subgroups of Lie groups}, e questo articolo ci suggerì di considerare nel disco iperbolico, su cui agisce $SL(2,\setR)$, un'orbita (simile ad un albero) del sottogruppo discreto $SL(2,\setZ)$, la cui struttura algebrica è quasi la stessa di un gruppo libero, e quindi essenzialmente di immergere un albero nel disco iperbolico. A questo punto, intendevamo introdurre su questo albero un'opportuna passeggiata casuale alla maniera di Furstenberg, che purtroppo procede in maniera non costruttiva, e di usarla per descrivere una famiglia di rappresentazioni del gruppo libero in analogia a quelle di $SL(2,\setR)$, tutte ben note e classificate. Si trattava di mescolare analisi armonica, geometria iperbolica e teoria della probabilità, una prospettiva stimolante di avanzamenti in ambiti molto di moda. In quel seminario cercammo di raggiungere questo obiettivo ricoprendo il disco iperbolico con i traslati di un dominio di Siegel sotto l'azione di $SL(2,\setR)$. Speravamo, per questa via, di ottenere, in modo esplicito e particolarmente regolare, alberi adatti a trovare una scelta semplice, quindi matematicamente trattabile, di un'opportuna passeggiata casuale, ossia di un opportuno operatore di transizione che mimasse l'operatore di Laplace--Beltrami sul disco iperbolico. Ma non ci riuscimmo, ed in effetti non era quella la via. Però questo mi rammentò la bella conferenza che Pierre Eymard aveva presentato al congresso organizzato da Sandro all'INdAM pochi anni prima, nel 1976, dove spiegava in maniera chiara ed avvincente risultati profondi di Sigurður Helgason sulle autofunzioni dell'operatore di Laplace--Beltrami sul disco iperbolico in termini del nucleo di Poisson ad esso associato. Quella fu la prima volta che dai contenuti di quella conferenza e dal bell'articolo espositivo che li presentava ricevetti un'illuminazione (la seconda fu più di quaranta anni dopo, in collaborazione con Wolfgang Woess e Maura Salvatori, ma in un contesto del tutto diverso). Fu questo che, traslato dal continuo al discreto, ci diede l'idea d'impiegare il nucleo di Poisson dell'operatore di transizione isotropo di passo uno su un albero omogeneo (che talvolta si chiama operatore di Laplace sull'albero) per introdurre rappresentazioni del gruppo libero a partire dalla sua frontiera di Poisson. Il nucleo di Poisson su alberi, non solo omogenei, era stato studiato dal punto di vista combinatorio in alcuni profondi lavori di Pierre Cartier. Quest'idea ci fornì la chiave per la costruzione di una famiglia analitica di rappresentazioni unitarie o uniformemente limitate del gruppo libero (le rappresentazioni sferiche) con proprietà simili a quelle delle omonime rappresentazioni di una famiglia diversissima: il gruppo $SL(2,\setR)$ e i gruppi di Lie semisemplici di rango uno.

Una volta stabilita la linea su cui procedere, Sandro e io dovemmo studiare molti altri settori, seguendo non solo l'articolo di Eymard ma anche i suggerimenti di nostri amici matematici profondi come Jacques Faraut ed in seguito Michael Cowling e Ádám Korányi: ad esempio le coppie di Gel\ms{}fand e le funzioni sferiche. Fu un periodo ``matto e disperatissimo'', un forsennato affanno quotidiano, ma elettrizzante: riuscivamo giorno per giorno come per magia a estendere definizioni e proprietà prima considerate in tutt'altro ambito ad un gruppo discreto difficile da studiare. La nostra teoria delle rappresentazioni sferiche del gruppo libero si affiancò allo studio di altre famiglie di rappresentazioni di quel gruppo considerate in quegli stessi anni, ad esempio da Tadeusz Pytlik e Ryszard Szwarc, ma con una base geometrica più intrinseca, che prometteva ulteriori significativi sviluppi, ad esempio generalizzazioni dal gruppo libero all'intero gruppo degli automorfismi di un albero omogeneo (sette anni dopo Ryszard parlò a lungo di matematica con Sandro in un anno che passò alla Sapienza come professore visitatore; Tadeusz in seguito visitò me a Tor Vergata e scrivemmo un articolo sulla disuguaglianza di Haagerup; entrambi collaborarono con Anna~Maria ed Anna sulle rappresentazioni del gruppo libero). Sandro ed io pubblicammo nel 1982 il nostro articolo sulle funzioni sferiche sul gruppo libero \cite{33}, nel quale dimostrammo l'irriducibilità delle rappresentazioni sferiche. Ampliammo poi considerevolmente quei risultati pubblicando l'anno successivo un libro \cite{34} che diede un forte impulso alla teoria delle rappresentazioni di questo gruppo assai ostico, spiegò per questa via lo spettro in $L^p$ del suo Laplaciano isotropo ed ispirò molti significativi sviluppi; vi spiegammo anche la magistrale analisi degli spazi di Lorentz e gli operatori di intrallacciamento con la quale Anna~Maria ed Anna completarono il nostro lavoro seguendo i consigli e le vaste conoscenze di Michael. Io ero un po' preoccupato della responsabilità di scrivere un libro, mentre Sandro ne era entusiasta: era sempre stato convinto che l'eredità scientifica di un matematico si basa sulle ampie trattazioni possibili in un libro più che su singoli articoli anche se molto influenti. In effetti, scrivere un libro porta ad una visione più ampia e globale di quella necessaria per un breve articolo scientifico, ma soprattutto spinge a riorganizzare la presentazione in maniera più organica e significativa. Nella stesura del testo fianco a fianco, imparai da Sandro la sua capacità inestimabile di presentare la matematica in modo che le idee emergano chiaramente dai pur indispensabili aspetti tecnici, e questo cambiò da quel momento in poi il mio modo di scrivere sia libri sia articoli. In effetti, Sandro aveva una capacità meravigliosa di comunicare le idee matematiche, ed il suo modo di scriverle rivelava non solo quanta gioia gli desse sviluppare quelle idee ma anche condividerle: questo lo rendeva in grado di esporle in maniera avvincente. A questo proposito, si vedano vari suoi articoli parzialmente espositivi \cites{17,25,
31,36,37,41,42,44,46,49}, oltre naturalmente alle monografie \cites{34,39,40}. Il fatto che frequentemente Sandro includesse sezioni espositive nei suoi articoli di ricerca è una testimonianza della sua passione per la comunicazione e condivisione della matematica.

I temi del libro, poi, spinsero Sandro e me a interessarci ad edifici (\emph{buildings}) di Bruhat--Tits ed ai loro legami con $SL(2,\setZ)$, brillantemente spiegati nel libro di Jean-Pierre Serre, \textit{Arbres, amalgames, $SL_2$} (anche qui, Sandro traduceva in italiano \emph{buildings}, o \emph{immeubles}, con ``palazzi'', io con ``edifici'', e non raggiungemmo mai un accordo nonostante varie disquisizioni filologiche). Sfruttando il nostro teorema di irriducibilità delle rappresentazioni sferiche e immergendo nel gruppo $p$-adico non discreto $PSL(2,\setQ_p)$ un sottogruppo discreto che agisce in maniera semplicemente transitiva su un albero, Sandro e io dimostrammo nel 1984 un risultato per noi sorprendente, ossia che ogni rappresentazione sferica unitaria di $PSL(2,\setQ_p)$ si restringe in maniera irriducibile a questo sottogruppo discreto \cite{35}. Vari risultati in ambiti connessi furono ottenuti da diversi altri studiosi nei decenni successivi, in particolare da Tim Steger, come vedremo fra poco.

In quel periodo, Sandro diresse le tesi di laurea di due studenti molto bravi, com i quali, come con tutti gli altri, scambiò un profondo affetto: Roberto Scaramuzzi ed in seguito Beatrice Pelloni. Entrambi conseguirono il dottorato a Yale; poi Beatrice fu in Grecia (IACM-FORTH, Heraklion), poi all'Imperial College, a Reading ed infine a Heriot--Watt University ad Edinburgo. E voglio citare due giovani brillanti che furono studenti miei, non di Sandro, ma la cui relazione di amicizia e la capacità profonda di scambio affettivo fu ugualmente importante anche con lui: Alessandra Iozzi (coetanea di Roberto, poi studentessa di dottorato a Chicago con Robert Zimmer) e Fausto Di Biase (coetaneo di Beatrice, poi studente a Saint~Louis con Steven Krantz). Voglio anche ricordare Martin Moskowitz, professore al Graduate Center della City University di New York (CUNY), che venne tante volte a visitarmi a Roma ed intrattenne una viva amicizia con Sandro e me. E rammento anche un giovane che era stato uno studente molto bravo alle lezioni di Analisi Matematica che svolgevo a Fisica, Paolo Emilio Barbano, il quale si trasferì al corso di laurea in Matematica, dove Sandro ed io congiuntamente, con grande pazienza e grande affetto, lo aiutammo a redigere la tesi ed a maturare senso critico e poi, grazie all'aiuto di Martin, lo avviammo agli studi dottorali a CUNY. Cito anche Alessandra Gallinari e Paolo Piccione, che ho avuto entrambi come brillanti studenti, esattamente come Fausto, ad un corso abbastanza avanzato in analisi matematica al Dipartimento ``Castelnuovo'', nel quale sotto mentite spoglie insegnavo anche un po' di analisi di Fourier classica; poi Alessandra e Paolo scrissero la tesi di laurea con Sandro, e dopo entrambi continuarono all'estero gli studi post-laurea: Alessandra conseguì il dottorato sotto la guida di Jonathan Rosenberg nella stessa università dove avevo ottenuto io stesso il Ph.D. e poi andò all'Università di Pennsylvania e dopo a Madrid, e Paolo divenne un'autorità nel mondo matematico brasiliano. E poi Fabio Rossi e Daniela Stalteri, che si laurearono con Sandro e proseguirono gli studi rispettivamente a Chicago e a Saint~Louis.

Quello del 1984 fu il nostro ultimo articolo insieme. Infatti, da quel momento l'attenzione di Sandro, sempre fedele all'ambito dell'analisi armonica, si rivolse esclusivamente alle funzioni sferiche per opportuni operatori di transizione su alberi omogenei ed alla conseguente teoria delle rappresentazioni del loro gruppo di invarianza, il gruppo libero oppure un analogo prodotto libero, ed alle passeggiate casuali, ovvero ai processi di diffusione, indotti da operatori di transizione sui vertici dell'albero ed invarianti rispetto a questi gruppi. Invece i miei interessi si spostarono gradualmente verso ambiti con proprietà di simmetria e di invarianza molto minori, ad esempio alberi non più omogenei od operatori di transizione non più invarianti, e quindi divennero distanti da quelli di Sandro, il quale fu rattristato da questa mia diaspora, e questo dispiacere gli rimase vivo anche in seguito (il che mi addolorò allora e dopo). Forse il cambiamento di ambito, da quello invariante sotto l'azione di gruppi a quello generale, merita alcuni commenti storici, perché diede luogo a nuove linee di ricerca non solo mie ma anche, simultaneamente, di una parte del gruppo di analisi armonica. Provo a delinearne gli inizi in base ai ricordi della mia evoluzione scientifica.

Abbiamo osservato che gli articoli pionieristici di Cartier riguardano alberi e operatori molto generali, senza simmetrie. Io personalmente fui stimolato da questo e dalla collaborazione con (e dalle idee di) un brillante giovane matematico austriaco che ho già nominato, Wolfgang Woess, il quale nel 1984--85 visitò come post-doc l'Università di Roma ``La Sapienza''. Rapidamente, Wolfgang divenne uno dei miei amici più cari e dei miei collaboratori più abituali ed in seguito fu Professore Associato e poi Ordinario a Milano (dove diede un notevole impulso all'approccio probabilistico e di teoria del potenziale all'analisi armonica su alberi e grafi rispetto a operatori di transizione via via più generali e sempre meno invarianti). Fui allo stesso tempo stimolato a questi ambiti senza simmetrie, ossia su alberi arbitrari, non omogenei, o anche grafi ed alcuni domini frattali, dalla interazione con altri due matematici più maturi con i quali sviluppai una profonda amicizia ed una continua collaborazione: Mitchell Taibleson ed Ádám Korányi alla Washington University (poi Ádám si trasferì al Lehman College del CUNY), che ho già nominato più sopra. Vennero più volte in visita a Roma e furono amici anche di Sandro; di più, Ádám fu spesso nostro commensale, e da persona di elevata cultura apprezzò sempre moltissimo quelle cene con Sandro in cui il piacere delle pietanze si univa a quello delle idee.
Nella transizione dall'analisi armonica basata su gruppi di simmetrie all'analisi in ambiti senza simmetrie, ossia la teoria del potenziale e la probabilità, Wolfgang ed io fummo ispirati dai lavori di vari matematici, inizialmente il suo maestro Peter Gerl, il succitato Cartier e poi Yves Derriennic; ma qui voglio citarne esplicitamente un altro, Vadim Kaimanovich, un giovane matematico russo di grande profondità che avevamo conosciuto nel 1989, solo pochi anni dopo che, poco più che ventenne, era diventato una autorità mondiale in teoria del potenziale, e precisamente nel settore delle frontiere di passeggiate casuali su gruppi discreti. Poi Vadim, quando Wolfgang si trasferì a Milano ed in seguito a Graz, collaborò a lungo con lui e fu a Milano come professore visitatore --- anzi quella, nel 1990, fu proprio la prima volta che uscì dall'Unione Sovietica. A Milano Wolfgang e Vadim scrissero il loro primo articolo congiunto sulle frontiere di passeggiate casuali, che fu l'inizio di una collaborazione costante. Vadim ha continuato a collaborare con analisti armonici italiani, se non come coautore almeno tramite scambi di idee, ad esempio un po' con me, ed in anni recenti intensamente con Fausto. Nel 1992 Vadim, che con gli analisti armonici italiani si sentiva in grande sintonia, vinse in Italia un concorso a Professore Associato, ma poiché era cittadino sovietico il Ministero della Pubblica Istruzione, con una delibera presa dopo tempi lunghissimi, decise che veniva meno il requisito giuridico della reciprocità ed annullò la sua vittoria concorsuale. In quella circostanza, Sandro, Leo ed io spendemmo molte energie per far avere a Vadim la sua posizione di professore: ecco un altro esempio del costante impegno di Sandro per rendere l'ambiente universitario italiano più internazionale; in questo caso, si aveva la percezione che la delibera ministeriale fosse stata forse influenzata da pressioni campaniliste di gruppi accademici che miravano a favorire i propri candidati rispetto ad altri più bravi ma stranieri. Sandro, anche prima di questi accadimenti, osservava sempre che era stato controproducente da parte del nostro governo limitare le acquisizioni di professori stranieri di valore in base al fatto che i loro governi non avrebbero fatto la stessa cosa con i nostri scienziati.

A proposito di Wolfgang, che fu anch'egli amico di Sandro, mi piace ricordare questo aneddoto: nell'anno in cui fu a Roma, lui, Claudio ed io andavamo quotidianamente a pranzo con Sandro, in una specie di ristorante vicino all'Università (veramente una ``panineria'', o ``paninoteca'', ma esito a scrivere questi neologismi), il quale preparava buoni panini ma ci metteva molto tempo. Wolfgang, da bravo austriaco ci teneva molto alla puntualità, soprattutto per quanto concerne l'ora dei pasti, era sempre un po' a disagio per l'attesa. I tempi poi si allungavano a causa di quotidiane incombenze impreviste ed urgenti di Sandro o mie che saltavano fuori all'ultimo minuto. Quindi ci avviavamo verso il ristorante quasi sempre con considerevole ritardo e qualche disappunto di Wolfgang; ma condividere i fastidi con gli amici li attenuava, ed almeno l'attesa per la preparazione del pranzo ci regalava piacevoli conversazioni. In quegli anni e nel decennio successivo, da Yale venne varie volte in visita Persi Diaconis, che era stato inizialmente un prestigiatore e dopo era diventato un matematico eminente in probabilità e statistica: ricordo le sue conferenze su quante volte fosse necessario tagliare un mazzo di carte per arrivare ad una distribuzione pienamente casuale. Persi applicava metodi di teoria delle rappresentazioni dei gruppi all'analisi statistica, come accenneremo anche dopo. Sandro aveva con lui grande affinità ed affiatamento.

Un'altra ragione per cui Sandro e io non scrivemmo altri lavori congiunti risiede nel fatto che io lasciai la Sapienza: divenni Professore Ordinario e mi trasferii prima all'Università dell'Aquila e poi a Roma ``Tor~Vergata''. Passai anche alcuni anni negli USA come Professore Visitatore, così come lui, ma non negli stessi periodi: quindi ci vedemmo di meno. Durante il mio triennio aquilano, al ritorno da quella sede visitavo Sandro una volta alla settimana al suo Dipartimento, che era vicino alla fermata dove il mio autobus dall'Aquila arrivava nel primo pomeriggio, quindi continuammo a parlare di matematica, oltre che di politica universitaria, e lo aiutai un pochino nella direzione della tesi di alcuni suoi studenti. In effetti, una sua studentessa di tesi, Laura Atanasi, continuò poi gli studi di dottorato sotto la guida mia e del mio stretto collaboratore Enrico, che si era appena trasferito alla Sapienza, ed in seguito lei scrisse alcuni articoli con me; lo stesso successe ad un'altra studentessa di Sandro, Federica Andreano, che portò a termine il suo dottorato in USA e dopo ritornò a Roma, anche a Tor~Vergata dove scrisse articoli con me. Negli anni successivi, Sandro guidò numerosi altri studenti brillanti, come Fabio Scarabotti, Filippo Tolli e Tullio Ceccherini-Silberstein. Il primo, di cui Sandro aveva una grande ammirazione, rimase a Roma, gli altri due furono da lui avviati a conseguire dottorati in prestigiose università statunitensi; tutti e tre insieme, spronati dall'interesse di Sandro, furono poi coautori di varie importanti monografie scientifiche. Anche con tutti loro ebbe un rapporto di affetto e stima che durò tutta la vita. Sandro aveva un senso fortissimo della indipendenza e della libertà di ricerca ed evitava di fossilizzarne gli ambiti: è significativo che spingesse sempre i suoi migliori allievi a continuare gli studi all'estero, ma mai a farlo necessariamente in sedi o sotto la guida di matematici orientati sui suoi stessi temi.

Dopo il 1990 mi trasferii a Tor~Vergata, e da allora la mia ricerca matematica si ridusse molto per mancanza di tempo. Ero occupatissimo come Direttore di Dipartimento, Senatore Accademico per circa 13 anni e creatore di corsi di laurea innovativi, ed ero sommerso dai contatti scientifici con le industrie e da una mole smisurata di attività didattica per tanti insegnamenti matematici esotici ed almeno altrettanti non matematici. Anzi, Sandro, che ha sempre dedicato un grandissimo interesse alla didattica della matematica, ed aveva un forte senso dello spirito, scherzava spesso su queste mie variazioni didattiche poco matematiche. Nello stesso periodo lui era impegnatissimo con la sua attività di organizzazione della politica scientifica, ed anche con periodi passati come professore visitatore alla Università di New South Wales, dove era Michael.

Fu quindi fisiologico che Sandro e io non continuassimo a pubblicare insieme, sebbene continuassimo a parlare di matematica insieme quando possibile. Questo però non fu dannoso alle ricerche sull'analisi armonica su alberi omogenei, perché in quel campo c'era ormai meno bisogno che scrivessimo articoli insieme: infatti lui ebbe la fortuna di trovare, dovrei forse dire quasi far nascere, un collaboratore molto più profondo di me nel condurre studi sulle rappresentazioni di gruppi liberi (e di molto altro): Tim Steger. Sandro lo incontrò nel 1983, durante un semestre passato come professore visitatore al Dipartimento di Matematica di Washington University, dove Tim era uno studente di dottorato che tutti ammiravano come brillantissimo ma che non aveva ancora trovato la sua via più congeniale verso il completamento degli studi. In effetti, quell'intero Dipartimento, preoccupato che Tim potesse abbandonare gli studi, fu entusiasta della sintonia e della stima reciproca che si svilupparono immediatamente fra Tim e Sandro e che durò poi per sempre, con una profonda e calorosa amicizia. Sandro raccontava una cosa che lo colpì subito di Tim: il modo in cui lui seguiva le lezioni. In esse, Sandro proponeva agli studenti problemi la cui risposta non era nota, diciamo piccole congetture; ogni volta Tim, il giorno dopo, si presentava con la dimostrazione oppure con un controesempio. D'altra parte, di primo acchito Tim ebbe di Sandro questa impressione, che mi fu raccontata in seguito da Fausto e del quale riporto testualmente le parole: ``l’incontro con Sandro mandò in frantumi lo stereotipo che lui, come molti statunitensi, aveva degli italiani, perché Sandro si presentò subito come persona coltissima, dall’ottimo inglese, sicuro di sé, signorile e con qualcosa di aristocratico nei modi''. Chiunque abbia conosciuto Sandro non può che condividere totalmente questa impressione.

Tim cominciò a scrivere la sua tesi di dottorato sotto la guida di Sandro: studiò le rappresentazioni sferiche su un gruppo libero munito di un operatore di transizione di passo uno, invariante per traslazione ma non necessariamente isotropo, divenne prestissimo autonomo ed ottenne risultati inaspettati e profondi, fra i quali una dimostrazione completa della algebricità del nucleo di Green, ed una elegante e potente generalizzazione al caso non isotropo dell'irriducibilità delle rappresentazioni sferiche che fornisce a questo teorema una dimostrazione non solo più generale, ma anche migliore di quella che nel caso isotropo avevamo dato Sandro ed io. Nel corso della loro collaborazione sulla sua tesi, Tim ha ricevuto moltissime idee da Sandro. Tim osserva che la maggioranza di esse risultavano sbagliate\dots\ ma che questo sottolinea uno dei principali punti di forza di Sandro: generava proficuamente idee e le esponeva senza nessuna paura, senza nessuna timidezza. Fra le tante idee che erano valide, due sono state la chiave per pressappoco tutto il lavoro di Tim negli anni successivi, tutte e due già menzionate sopra: (1) che nello studio del gruppo libero o di un albero, il ruolo della loro frontiera è fondamentale; (2) che vale la pena capire bene i \emph{buildings} affini, di cui l'albero omogeneo rappresenta il caso più semplice.

Infine, i contenuti della tesi di Tim furono redatti e pubblicati congiuntamente con Sandro in un memorabile volume dei Memoirs of the American Mathematical Society \cite{40} (il gioco di parole è intenzionale) che apparve a stampa purtroppo più di dieci anni dopo l'inizio di quel lavoro di tesi, il quale, più che una tesi, si rivelava già come un lavoro di piena maturità. Durante la lunga gestazione Tim spiegò i suoi risultati in una serie di splendidi seminari a Roma ``La Sapienza'', le cui note da lui succintamente redatte rimasero per anni l'unico testo disponibile per questa teoria.

Dopo il dottorato a Saint Louis, Tim passò tre anni a Yale, quattro a Chicago, due all'Università della Georgia. Poi si stabilì come Professore a Sassari, e per varie ragioni qualche anno dopo ridusse i propri spostamenti, ma continuò a scrivere articoli di grande impatto in molte aree, e certamente, in particolare, sulle rappresentazioni di gruppi liberi, tematica nella quale è un indiscusso esperto insieme a Gabriella Kuhn; altri che hanno collaborato su quella tematica sono stati Waldemar Hebisch, Alessandra Iozzi e Sandra Saliani (Waldemar non collaborò mai direttamente con Sandro, ma si dimostrò tante volte disponibile ed utile a diversi altri componenti del gruppo). Tim è anche un indiscusso esperto su gruppi che agiscono su \emph{buildings} affini (suoi frequenti collaboratori sono stati Donald Cartwright, Anna Maria ed Anna) ed ha ottenuto importanti risultati sulla rigidità di restrizioni a reticoli di rappresentazioni di $PSL(2,\setR)$ (collaborando con Christopher Bishop) ed estendendo poi questo risultato a tutti i gruppi di Lie semisemplici (con Michael). A Sassari collaborarono con Tim o discussero temi di matematica anche Carlo Andrea Pensavalle, che conseguì la laurea magistrale a Genova nel 1986 sotto la direzione di Anna~Maria e poi raggiunse a Sassari Tim, con il quale scrisse vari articoli interessanti sulla irriducibilità della restrizione ad un sottogruppo libero di certe rappresentazioni del gruppo degli automorfismi di un albero omogeneo, ed in seguito su famiglie di rappresentazioni del gruppo libero. La collaborazione con Tim, penso più di qualsiasi altra, fu per Sandro fonte di grande gioia e di orgoglio per il suo ruolo sia di maestro sia di matematico.

Sandro continuò a occuparsi di rappresentazioni di gruppi di automorfismi di alberi omogenei insieme a Claudio Nebbia, studiando in questo contesto la teoria di Grigoriĭ I. Ol\ms{}shanskiĭ delle rappresentazioni cuspidali, la quale estende a gruppi discreti semplicemente transitivi su alberi, come il gruppo libero, l'analogia con le varie classi di rappresentazioni di $SL(2,\setR)$, non solo quelle sferiche (ossia quelle che ammettono un vettore $K$-invariante: la serie principale e complementare) ma anche le altre, ad esempio la serie discreta. Come già accennato, Sandro e Claudio scrissero su questi argomenti un altro libro, che risultò di grande influenza sulla letteratura matematica \cite{39}. Claudio continuò a parlare di matematica con Sandro e subito dopo trovò interessanti risultati su gruppi di automorfismi di alberi omogenei legati all'amenabilità ed alla proprietà di Kunze--Stein. In quegli anni Sandro ebbe altri studenti di ricerca, come Flavio Angelini, che si laureò nel 1989, poi conseguì il Ph.D. a UCLA su argomenti di geometria algebrica (una versione algebrica della dimostrazione di Demailly delle disuguaglianze asintotiche di Morse su spazi di Kähler), poi andò a Nizza dove si interessò a temi di crittografia, ed infine spostò i propri interessi scientifici verso problemi di matematica finanziaria ed economia, che ora sviluppa come professore a Perugia.

Stiamo ormai parlando degli anni successivi al 1993, quando cominciai la mia attività, impegnativa fin quasi allo sfinimento, di direzione del Dipartimento di Matematica a Tor~Vergata e di organizzazione di quell'Università e della sua amministrazione e didattica: passavo in ufficio o in aula sei interi giorni alla settimana e non avevo più tempo di visitare la Sapienza tranne che per occasionali conferenze. Quindi di questo periodo non ho molti ricordi da presentare sull'attività di Sandro, che incontravo quasi solamente al Convegno annuale di Analisi Armonica; peraltro cominciò a collaborare con Enrico ed altri, fra i quali il mio collega a Tor~Vergata Paolo Baldi, ed ebbe ancora studenti, ad esempio Giovanni Stegel che si trasferì a Sassari a collaborare con Tim. Sandro scrisse anche alcuni articoli su processi di diffusione e variabili aleatorie stabili su campi locali e sul semipiano iperbolico congiuntamente con Paolo ed Enrico e Marc Yor \cites{43,45} e più in generale su spazi ultrametrici \cites{41,44}, ed anche con Mauro Del Muto \cites{46,47}, nonché un articolo di statistica \cite{49} in collaborazione con Alberto Guerriero, Alberto Leone, Gian Piero Mignoli ed Enrico Rogora; con quest'ultimo coautore scrisse anche un libro didattico per studenti su serie di Fourier e di Fourier--Walsh.

A questo proposito, è significativa la genesi del suo articolo di statistica e di questo quaderno di note per studenti, entrambi scritti con Enrico Rogora, che cito quasi letteralmente --- le sue parole illustrano la dedizione e l'innovatività di Sandro verso la didattica. Alla fine del decennio di fine secolo Sandro aveva deciso di insegnare i corsi di calcolo al corso di laurea in Informatica. C'era una grossa sproporzione tra il numero degli studenti che frequentavano il corso ed il numero di coloro che ne superavano gli esami. Questo fatto incoraggiava a rivedere i programmi e le modalità d'esame. Enrico Rogora aveva sviluppato un programma per computer che a partire da una banca dati di domande a risposta multipla produceva versioni personalizzate degli scritti cambiando le domande scelte per ogni sezione e rimescolando le varie risposte di ciascuna domanda: oggigiorno troviamo facilmente in commercio versioni commerciali di tali applicativi, ma a quell'epoca non esistevano ancora. In questa maniera varie centinaia di studenti si potevano affollare nelle aule d'esame senza che fosse facile copiare le risposte.
Enrico Rogora discuteva con Sandro le statistiche elaborate dai dati raccolti. Durante quelle discussioni Sandro gli propose di mettere alla prova un metodo di indagine statistica che gli era venuto in mente leggendo il libro di Persi Diaconis sui metodi della rappresentazioni dei gruppi nell'analisi dei dati. I dati raccolti con gli esami non erano l'ideale per questa analisi, e quindi Sandro prese contatti con Alma Laurea per ottenere un set di dati più significativo. Entrambi andarono ripetutamente a Bologna ad incontrare il gruppo di Alma Laurea, che in effetti fu poi coautore dell'articolo~\cite{49}, con cui scelsero i dati da elaborare, prepararono gli script per implementare l'algoritmo suggerito da Sandro e discussero i risultati.

Anche il quaderno sull'analisi di Fourier--Walsh nasce dall'esperienza di insegnamento ad Informatica. Esso illustra un aspetto caratteristico del suo insegnamento, che intendeva guidare l' apprendimento degli studenti attraverso l'uso di una serie di domande a risposta multipla che focalizzassero i punti principali, stimolassero la riflessione su alcuni problemi scelti e facilitassero la comprensione degli argomenti e degli esercizi proponendo soluzioni dettagliate delle risposte esatte e il commento degli errori principali che commettevano.

In effetti, Sandro ha sempre avuto un grandissimo interesse ed entusiasmo per la didattica della matematica, ed è stato un eccellente docente. Dopo che a fine 2010 concluse il suo lunghissimo periodo come professore universitario, tenne ancora insegnamenti a contratto per tre anni nelle tre Università statali di Roma, fra cui un corso di Analisi Armonica che avevo disegnato a Tor~Vergata per far acquisire agli studenti non matematici gli strumenti utili per le applicazioni industriali: dopo tanti anni, quella fu la prima volta che Sandro insegnò analisi armonica mirata alle applicazioni e fui molto contento di condividere con lui quegli aspetti innovativi della nostra disciplina.

Sandro ha sempre anche avuto una passione innata per la ricerca. Aborriva la ricerca intesa come la produzione di minuscole estensioni di risultati già noti sotto ipotesi lievemente più generali (e di solito più tecniche e meno naturali), che designava come provinciale. Per lui la ricerca era l'esplorazione di un universo mentale sconosciuto, che ogni volta occorreva esplorare e scoprire, guidati dall'idea, ossia dal problema. Infatti, trovava una ideale corrispondenza fra la matematica e la filosofia di Platone di un mondo di idee (anche qui, l'uso dell'aggettivo ``ideale'' è un gioco di parole intenzionale: il duplice impiego si riferisce a significati diversi, ma affini). Proprio per questo aspetto del proprio pensiero, asseriva sempre, come già detto, di non sapere nulla di matematica e di affrontare ogni volta un terreno per lui incolto (ma proprio per questo affascinante). A questo proposito, includo una sua citazione autografa che riguarda questo punto. Si tratta di un brano di una conversazione che ebbe con Fausto nel 2009, a proposito di un risultato inaspettato che Fausto aveva scoperto tre anni prima. Il risultato era l'indecidibilità nel senso di Gödel di un problema aperto esattamente da un secolo: il problema di Fatou. Questo problema consisteva nel determinare se esiste sempre o no, senza ipotesi particolari, il limite asintotico di funzioni armoniche sul disco unitario quando si raggiunge la frontiera lungo percorsi tangenziali. Sotto vari casi diversi di ipotesi restrittive, nel corso di quei cento anni si era dimostrata l'esistenza del limite, caso per caso; sotto altre ipotesi restrittive, si erano trovati controesempi; ma in generale la risposta era ignota. Fausto, con l'acume e la tenacia che lo contraddistinguono, motivò a questo problema un team di coautori, alcuni dei quali esperti in teoria del potenziale ed altri in logica, ed insieme dimostrarono che il problema era indecidibile: la risposta dipendeva dalla scelta dell'insieme di assiomi che si sceglie di adottare per descrivere le basi concettuali della matematica, normalmente gli assiomi di Zermelo--Fraenkel e Assioma della Scelta (ZFC; qui C sta per \emph{Choice}). Fu, secondo Sandro ed anche me, un risultato epocale, che non ebbe, a parere di Sandro ed assolutamente anche mio, la giusta rilevanza che meritava --- ed in effetti questo è un altro esempio in cui Sandro criticava il provincialismo campanilista di alcune parti dell'ambito matematico italiano, dove pochi studiosi di analisi matematica capiscono a fondo la logica matematica e ne percepiscono l'importanza. Ebbene, cito un brano del messaggio che Sandro inviò a Fausto riguardo a quel lavoro. Sandro scrisse:
``Stavo descrivendo il sentimento prevalente di chi fa ricerca matematica (o almeno il mio sentimento prevalente) che è quello di scoprire paesaggi sconosciuti e tentare di comunicare le scoperte ad altri. In effetti sono più che paesaggi, perché sono paesaggi popolati da entità più o meno mostruose, insomma avventure. Nel libro che abbiamo scritto con Massimo Picardello sull'analisi armonica su gruppi liberi, prima della prefazione citiamo una frase di Bernal Díaz~del~Castillo che egli premette alla sua `Historia verdadera de la conquista de Nueva~España':
\emph{Mas lo que yo oí y me hallé en ello peleando, como buen testigo de vista, yo lo escrebiré, con el ayuda de Dios, muy llanamente, sin torcer a una parte ni a otra}. Nella traduzione in inglese del libro di Diaz del Castillo (io ho letto il libro solo in inglese) questa frase è così tradotta:
\emph{That which I have myself seen and the fighting I~have gone through, with the help of God, I will describe quite simply, as a fair eye witness without twisting events one way or another.}
La matematica come racconto di un'avventura. Come dico ai miei studenti, un problema matematico da risolvere è come una lanterna che si può usare per illuminare la strada. Allora, da questo punto di vista (e sto parlando di emozioni non di pensiero razionale) se il problema si rivela non decidibile, la lanterna finisce per brillare ancora di più, anzi permette di guardare la strada, o l'avventura, dall'alto, in modo che si possano contemplare i lontani confini del paesaggio, che forse sono gli assiomi ZFC.''

Sandro ha avuto altrettanta passione per l'organizzazione della ricerca e dell'insegnamento, ossia per la politica scientifica. Ecco un'altra citazione approssimativa di una sua frase, riferitami da Tim: ``Quasi tutti i matematici che si trovano coinvolti nella politica e nell'amministrazione si lamentano in continuazione che ciò interferisce con il loro vero lavoro: la ricerca. Non è il caso per me, e forse neanche così tanto per loro. Se la ricerca non va bene, o se ci si deprime nel considerare il talento esagerato di tanti altri, questi altri impegni aiutano a mantenere l'equilibrio psicologico.''

Sandro era insofferente alla presunzione ed alla insincerità, che sapeva riconoscere immediatamente ed avversava con determinazione.
Come abbiamo detto, incoraggiava i suoi studenti ad approfondire gli studi presso altre scuole, principalmente negli Stati Uniti, al fine di fargli ampliare i loro orizzonti scientifici e magari differenziarli dai suoi, ma si aspettava sempre che loro studiassero ed insegnassero con la giusta umiltà, ad esempio dedicando la massima cura ed il giusto entusiasmo ai corsi che insegnavano: si infastidiva se qualcuno dei suoi allievi si riferiva alla didattica universitaria come \emph{a necessary evil}, una espressione frequentemente pronunciata da giovani ricercatori statunitensi un po' presuntuosi. In effetti, Sandro percepiva come arroganti alcuni ``economisti che hanno studiato negli Stati Uniti e si ritengono per questo superiore a chiunque altro, economista o politico'' (sue parole testuali), e non aveva peraltro molta stima della teoria dell’equilibrio economico generale; tuttavia, la consapevolezza dell’importanza civile di discutere con coloro che hanno una opinione diversa dalla propria e lo sforzo di arricchire il dibattito culturale e politico dei contenuti che sapeva trasmettere e dei dubbi che era in grado di introdurre lo spingevano a dialogare con loro, come con altri. Sandro aveva una profonda e vasta cultura classica, ma era diffidente verso coloro che chiamava ``parolai'', capaci di sedurre con l’arte retorica ma privi di una solida base di conoscenze reali.

Nella sua permanenza negli Stati Uniti Sandro aveva maturato ideali di libertà e democrazia. Quando nel 1969 prese servizio come professore in Italia, all'Università di Genova, manifestava tratti tipici del comportamento dei professori americani, ad esempio la grande amichevolezza con gli studenti. Di nuovo in linea con gli ideali americani, si sentiva \emph{liberal}, concetto che cercò di tradurre agli studenti dicendo di essere ``liberale''. Gli studenti sgranarono gli occhi: erano gli anni della rivoluzione studentesca, e nessuno studente aveva mai visto all'Università un liberale\dots

Sandro aveva un forte senso della giustizia e dell'equità, ed una grande passione per l'impegno civile e per la giustizia sociale. I suoi articoli pubblicati in quotidiani e riviste (recentemente raccolti in un'apposita pubblicazione) rivolgevano l'attenzione ai gravi problemi del sistema giudiziario italiano e le loro conseguenze per i cittadini. Esempio di questi problemi è dato dalle procedure inquisitorie nel triste caso dell'omicidio della studentessa Marta Russo, assassinata nel 1997 mentre camminava nel campus della Sapienza. Le indagini balistiche indicarono come punto di origine dello sparo la finestra di una stanza usata da due dottorandi di filosofia, immediatamente arrestati con grande compiacimento e soddisfazione dell'opinione pubblica. I due presunti assassini rimasero in custodia cautelare per un lunghissimo periodo, e in seguito furono condannati. La sentenza di colpevolezza si basava su un'unica testimonianza, di un'impiegata dell'università che frequentava la stessa stanza degli accusati. Sandro ed Alberto Beretta~Anguissola pubblicarono un libro, \textit{La prenderemo per omicida} (Koinè Nuove Edizioni, 2001, con prefazione di Giovanni Sabbatucci), che contiene la versione integrale dell'interrogatorio, dalla quale Sandro e Beretta~Anguissola evincono che la testimone abbia subito pressioni e minacce, e deducono che ne se sia stata spinta a firmare una testimonianza falsa.

Altrettanto grande era la sua passione per l’impegno civile e per la giustizia sociale. Poco prima della propria morte, Carlo Pucci volle affidare a Sandro la cura della Fondazione da lui creata e intitolata congiuntamente a Ernesto Rossi (suo zio da parte di madre) ed a Gaetano Salvemini. Tale eredità comprendeva la tutela dei due importanti archivi personali di Rossi e Salvemini, il primo depositato presso gli Archivi Storici dell’Unione Europea a Firenze e il secondo custodito dall’Istituto Storico della Resistenza in Toscana ma sino ad allora poco valorizzato. Sandro non risparmiò energie per risolvere i problemi burocratici e pratici incontrati nell’eseguire il compito affidatogli. Grazie al suo impegno riuscì a sistemare le questioni relative alla custodia e alla valorizzazione degli archivi della Fondazione e presiedette con onore il suo Consiglio direttivo dal 2003 al 2016. Negli anni della sua presidenza, la Fondazione Ernesto Rossi--Gaetano Salvemini, grazie all’entusiasmo e al lavoro volontario dei suoi componenti, ha difeso e diffuso gli ideali di democrazia, giustizia e libertà che furono propri dei suoi maggiori, attraverso la presentazione di libri, convegni, dibattiti e con premi a tesi di laurea e dottorato attinenti a questi temi. Grazie al suo esempio e nel suo ricordo, l’attività della Fondazione prosegue ancora oggi lungo le direttive impostate da Sandro e da Pucci.

Sempre in tema di giustizia ed equità, un altro tratto di Sandro forse rafforzato dalla esperienza americana era il suo forte senso delle garanzie personali (il neologismo corrente è ``garantismo''). A partire dal 1992, il parlamento e la magistratura italiani approvarono lo strumento della ``custodia cautelare'', che permetteva alla polizia di tenere in carcere persone sulle quali aveva sospetti, forse anche come strumento di pressione psicologica per estorcergli confessioni. Sandro protestava vivamente contro questa pratica, ed osservava che, se questo fosse avvenuto negli Stati Uniti, nell'arco di poche ore qualche giudice avrebbe emesso una ingiunzione di \emph{habeas corpus}, liberando immediatamente il carcerato la cui detenzione era stata inflitta senza concreti elementi di accusa. Questa pratica fu utilizzata insistentemente, con grande supporto popolare, a partire dal 1992, ovvero dall'inchiesta ``Mani Pulite'' presso la Procura della Repubblica di Milano, il cui Procuratore capo, Francesco Saverio Borrelli, sovrintendeva una squadra di procuratori aggiunti molto agguerriti e secondo alcuni un po' spregiudicati. Sandro aveva coniato una definizione ironica per i sostenitori di Borrelli: con un divertente gioco di parole purtroppo comprensibile solo dai matematici, li chiamava ``la tribù dei Borrelliani''.

A proposito dei punti di vista rinforzati dalla esperienza americana: vari tratti della attività organizzativa e didattica di Sandro si richiamavano a procedure usuali negli Stati Uniti ma ancora non in Italia. Ad esempio, l'utilizzazione frequente di test d'esame stringati, talvolta (negli ultimi anni) quiz a risposte multiple, e l'insistenza sul risolvere esercizi. In effetti, appena arrivato in Italia come professore, propose ed ottenne la pubblicazione presso l'editore Boringhieri di un libro di testo americano di successo sull'Analisi Matematica di base, il libro di Tom Apostol, che rovesciava in maniera intelligente l'ordine espositivo consueto in Italia e conteneva molti esercizi, che erano parte integrante, non accessoria, dell'insegnamento: ne curò anche l'edizione italiana. Molti anni dopo, come direttore di Dipartimento, portò avanti un altro interessante progetto
: quello di introdurre nei corridoi dell'edificio fontanelle di acqua potabile, comunissime ovunque negli Stati Uniti ma mai viste in Italia, ed in particolare negli edifici universitari romani, privi di condizionamento d'aria nei corridoi e nelle aule, ma torridi d'estate. Sandro si vantava molto di questa sua realizzazione (anche se la posa in opera fu poi compiuta dal suo successore alla guida del Dipartimento, Enzo Nesi). In effetti, egli fu un precursore: adesso le fontanelle nei corridoi universitari stanno diventando molto più comuni.

L'amichevolezza di Sandro verso tutti matematici del suo gruppo non significava che non ci fossero dissensi. Ce ne furono, ma soprattutto con altri professori, al di fuori del gruppo di analisi armonica, per cui Sandro non aveva simpatia. Tim Steger, che fu un grande amico di Sandro, direi una persona per lui di fondamentale importanza affettiva oltre che scientifica, mi ha recentemente riferito un giudizio che riporto (le parole, ricordate dopo tanti anni, non sono precise, ma il senso lo è):
``A un certo punto ho riconosciuto che alcuni matematici di grande talento erano anche persone sgradevoli e arroganti che trattavano male gli altri, dentro e fuori il mondo della matematica. Ciò ha notevolmente ridotto la mia ambizione di essere uno di questi grandi.''
Sandro era sincero e franco e non esitava ad esprimere le proprie opinioni, anche se antitetiche a quelle di altri. Soprattutto, non esitava nell'esprimere apertamente i suoi punti di vista scomodi su quali fossero i comportamenti opportuni in relazione al pubblico interesse, alla modernizzazione del mondo accademico italiano, alla necessità della sua internazionalizzazione, alla necessità che il reclutamento accademico avvenisse sulla base del merito scientifico e non politico o sindacale, al superamento dei campanilismi locali. Spesso li esprimeva con battute sferzanti, ed è facile capire come fra alcuni suoi colleghi si generassero nei suoi confronti malumori, a volte antagonismi ed irritazioni.

Abbiamo già osservato che la cultura di Sandro fosse vasta. In effetti, aveva una grande vitalità ed una gran mole di interessi culturali, che lo affascinavano. Leggeva molto, e come abbiamo visto scriveva anche molti articoli su giornali e riviste. La sua curiosità e constante ricerca di nuove conoscenze lo portarono ad appassionarsi ad altri campi del sapere, come la psicoanalisi. Dopo il doloroso lutto per la perdita del fratello maggiore Lorenzo, cercò sollievo sottoponendosi alla psicoterapia. Aveva approfondito tutte le teorie psicoanalitiche, e strinse amicizia con alcuni cultori di psicoanalisi. Fu invitato a svolgere lezioni di matematica per i membri della Società Romana di Psicoanalisi. In questo ambito presentò e discusse con loro la sua convinzione circa l'esistenza di una analogia tra il lavoro psicoanalitico ed i processi mentali ed emotivi che avvengono nel corso dell’attività di ricerca in matematica. Queste sue idee furono oggetto di lunghe discussioni con i suoi amici psicoanalisti, che poi pubblicarono un articolo su questi temi in una rivista di psicoanalisi.

In effetti, anche molto in seguito in lui affioravano legami fra psicoanalisi e matematica: diceva che la psicoanalisi gli serviva a ``pulire i canali del piacere'', e quindi forniva impulso e motivazione per la sua ricerca matematica. Come dirò fra poco, negli ultimi anni della sua vita Sandro faceva un sogno ricorrente consistente di connessioni di idee matematiche relative al primo problema che aveva affrontato e risolto --- nel sogno, il problema si ampliava e quelle idee lo risolvevano. Gli ho chiesto spesso più ragguagli per capire se nella semi-coscienza onirica il controllo del suo ego garantisse a quelle idee sostanza e rigore, ma da sveglio gli era difficile ritrovarne le fila, come sempre accade in questi casi.


Sandro visse con qualche malinconia la fase del pensionamento. Appena pensionato, cercò di disfarsi di tutta la sua biblioteca matematica regalandola a vari istituti. Poco dopo che Sandro raggiunse l'età di settanta anni, il Ministero dell'Università e della Ricerca Scientifica aveva preso una decisione assai dubbia dal punto di vista della buona gestione dei fondi pubblici: per ridurre il budget necessario per il Ministero (a danni di quello dell'Istituto nazionale di Previdenza Sociale che paga le pensioni), il ministro aveva anticipato l'età del pensionamento dei docenti universitari, che per Sandro passava dagli originali settantacinque anni a settantadue, ed anzi all'inizio sembrò che si accorciasse subito a settanta, cosa che in seguito avvenne. Sandro si sentiva defraudato della sua aspettativa di continuare la sua attività didattica fino al termine previsto originalmente, ed esibiva il proprio dissenso contro questo fato maligno non solo lamentandosene, ma anche con azioni dimostrative. Ad esempio, fuori del suo ufficio al Dipartimento ``Castelnuovo'' aveva disposto con più di un anno di anticipo grandi scatoloni di cartone destinati a contenere i suoi libri ed appunti da portar via al momento della fine del servizio: questi scatoloni non entravano nel suo ufficio, che era molto piccolo, e venivano lasciati vuoti ma in piena visibilità come segno di protesta\dots

Qualche anno dopo il 2010, ovvero dopo il periodo degli ultimi suoi articoli scientifici e del pensionamento, il decorso di una malattia neurologica ed una grave invalidità limitarono sempre più i suoi interessi matematici, ma mai la sua lucidità e la sua curiosità, e mai, in nessun modo, la sua capacità di scambiare affetto, l'aver condiviso il quale è per me più importante dell'intero processo di formazione scientifica che grazie a lui ho avuto. Oltre all'affetto ed all'intelligenza, rimase sempre intatta la sua umanità, così come quella di Irene, ai quali entrambi sono grato di avermi aiutato, nonostante la malattia neurologica di Sandro, ad affrontare una mia situazione personale dolorosa di natura affine alla sua. A poco a poco la voce si affievoliva e non era più stentorea, nel discorrere non gli venivano sempre in mente le parole e talvolta perdeva il filo, verso la fine aveva momenti di tristezza che mi instillavano una punta di pena rammentandomi la sua abituale radiosità. La memoria dei fatti recenti si indeboliva. A volte, parlando con me, si dimenticava quali fossero i risultati che avevamo ottenuto insieme, e credo che si dimenticasse analogamente il suo lavoro scientifico di quasi tutta la vita, l'analisi armonica su alberi; si rammentava però l'analisi funzionale della sua tesi di dottorato, come dirò. Ma molte memorie antiche erano ancora vive; in effetti, l'ultimo suo lavoro, scritto subito prima dell'insorgere dei sintomi neurologici, non fu matematico, ma fu comunque su un albero: il suo albero genealogico, che Irene decise di redigere, cosa che fece grazie all'aiuto di Sandro, basato soprattutto sulla sua prodigiosa memoria, rafforzata anche dalla riesplorazione di documenti, lettere antiche e ricordi atavici. Fu in grado di ricostruire il ramo paterno fino al diciottesimo secolo, mi pare fin oltre la confluenza dei cognomi Figà e Talamanca. Io non vidi l'albero genealogico che aveva redatto, e d'altra parte leggere un albero genealogico non dice molto a chi non conosce già tutti i membri della famiglia, ma me ne parlò nelle fasi avanzate della malattia. Raccontava la storia dei suoi antenati che aveva conosciuto, e quella dei suoi bisavoli e trisavoli come gliela avevano raccontata i genitori, e così via generazione dopo generazione, con l'abituale lucidità e acume, ne era molto fiero e gli brillavano gli occhi.

Sempre, fino alla fine, ad ogni visita, in ogni momento, si ricreava quel piacere di vederci, di parlarci, di capirci che era stato una costante del nostro rapporto di mezzo secolo; e lo stesso succedeva con gli altri amici che lo visitavano. Per tale motivo questa introduzione al suo lavoro scientifico, che temo di essere inadeguato ad illustrare nella sua completezza e nelle sue conseguenze sulla letteratura matematica, è di fatto diventata un ricordo della sua personalità e del suo fascino.

Mi si lasci chiudere queste pagine con un ultimo aneddoto, che risale agli anni vicini al 1980. A quell'epoca Sandro era ben conscio di aver dato contributi significativi alla ricerca matematica e di aver creato dal nulla una scuola, ma era anche consapevole del suo ruolo e della sua importanza via via crescente nella politica scientifica (descriveva sé stesso e il piccolo gruppo di colleghi che collaboravano con lui nella politica della matematica come ``quarantenni d'assalto''). Perciò presumeva che, al termine della sua vita, su di lui sarebbe stato scritto un necrologio, ma, in parte scherzosamente ed in parte seriamente, diceva che non voleva che esso gli fosse scritto da qualche collega invidioso e forse neppure ben addentro alle sue molteplici attività, o persino da qualche allievo non sufficientemente capace. Quindi, almeno per la parte scientifica, se lo sarebbe scritto da solo poco a poco. In effetti aveva cominciato a scrivere una descrizione sintetica ma accurata dei propri articoli scientifici quando si preparava ad affrontare il suo concorso a cattedra, poi l'aveva continuata, ma in seguito smise, ed ora non trovo le pagine che aveva scritto (che mi sarebbero state utili adesso per descriverne più accuratamente l'opera scientifica).
Così, qui, forse un poco anche per il timore di inverare la sua preoccupazione (che un suo allievo non particolarmente capace gli scrivesse il necrologio) ho evitato di fare un'analisi ragionata e scientifica dei suoi risultati e del loro impatto, compito in effetti troppo difficile. Ma la ragione vera è che ho trovato più naturale tracciare uno schizzo della sua personalità, non solo scientifica ma soprattutto affettiva, per come la ricordo negli oltre 50 anni di amicizia. Credo che negli ultimi anni della sua vita non gli interessasse più un necrologio scientifico; la malattia gli stava portando via i ricordi dei risultati che aveva ottenuto, tranne il primissimo, che riaffiorava nel suo mondo mentale in maniera onirica: sviluppi possibili di quel risultato gli si manifestavano frequentemente in sogni, che spesso mi raccontava. Ancor prima dell'insorgere della malattia od ai suoi inizi, quando organizzammo vari congressi celebrativi per il suo 50esimo compleanno, poi il 70esimo, poi il pensionamento, poi l'80esimo, diceva che gli stavamo facendo la cerimonia funebre da vivo e non la voleva; tuttavia partecipava alle conferenze e tempestava di domande i conferenzieri. Invece, all'affetto è sempre rimasto sensibile sino alla fine e lo ha sempre condiviso con gli amici: probabilmente quindi gli sarebbe piaciuto essere ricordato con questo genere di memorie. Spero che suscitino un po' di emozione in chi lo ha conosciuto come l'hanno suscitata in me mentre le rievocavo.

Ci tengo molto a ringraziare Irene di avermi spinto a scrivere questa memoria. So che non lo ha chiesto a me solo per la mia lunga conoscenza di Sandro, e non ha inteso chiedermi di scrivere solo una analisi scientifica un po' distaccata dei risultati matematici. Sono certo che pensasse a pagine che ci facessero condividere i ricordi, di Sandro e con Sandro, ciò che, in effetti, mi è stato naturale scrivere: memorie di momenti passati insieme, raccontati con il sorriso, non solo con malinconia.

Desidero esprimere un vivo ringraziamento ad Enrico Casadio~Tarabusi, amico e coautore sia di Sandro sia mio, per aver dato un contributo entusiastico ed assolutamente sostanziale alla raccolta degli scritti scientifici di Sandro di cui questa biografia era l'introduzione. Anche gli aiuti di Michael Cowling, Fausto Di~Biase, Tim Steger e Leonede de~Michele, anche loro tutti amici comuni, sono stati preziosi. Tutti e cinque, inoltre, hanno contribuito ad avviarmi verso una maggiore precisione linguistica e storica, una abilità che, accanto alla ricchezza di affetto, li contraddistingue. Ringrazio Gianni Mauceri e Fulvio Ricci che mi hanno dato utili ed accurate indicazioni storiche difficilmente reperibili, oltre che la condivisione di affetto che traspare in queste pagine. Sono anche grato a Fausto e a Tim per avermi aiutato a delineare tratti della personalità di Sandro, avermi dato la gioia di scoprirne o riscoprirne alcuni, e di
avermi fornito alcune citazioni autografe di Sandro. Tutti poi hanno avuto la pazienza di compiere un'opera benemerita che certo non potevo aspettarmi né tanto meno chiedere: la correzione dei miei molteplici errori tipografici. Ma so che la loro vera ragione non è stata solo l'abituale gentilezza verso di me, ma anche il desiderio di rendere più curato, più perfetto
questo affettuoso ricordo di Sandro.
%
%


\renewcommand\refname{Pubblicazioni scientifiche di Alessandro Figà-Talamanca}


\begin{bibdiv}
%
\begin{biblist}
%
\DefineName{baldi}          {Baldi,            Paolo}
\DefineName{casadiotarabusi}{Casadio~Tarabusi, Enrico}
\DefineName{cecchini}       {Cecchini,         Carlo}
\DefineName{cohen}          {Cohen,            Joel~M.}
\DefineName{curtis}         {Curtis,           Philip~C., Jr.}
\DefineName{delmuto}        {Del~Muto,         Mauro}
\DefineName{demichele}      {De~Michele,       Leonede}
\DefineName{franklin}       {Franklin,         Stanley~P.}
\DefineName{gaudry}         {Gaudry,           Garth~I.}
\DefineName{guerriero}      {Guerriero,        Angelo}
\DefineName{leone}          {Leone,            Alberto}
\DefineName{mauceri}        {Mauceri,          Giancarlo}
\DefineName{mignoli}        {Mignoli,          Gian~Piero}
\DefineName{nebbia}         {Nebbia,           Claudio}
\DefineName{picardello}     {Picardello,       Massimo~A.}
\DefineName{price}          {Price,            John~F.}
\DefineName{rider}          {Rider,            Daniel~G.}
\DefineName{rogora}         {Rogora,           Enrico}
\DefineName{steger}         {Steger,           Tim}
\DefineName{yor}            {Yor,              Marc}
%
\DefinePublisher{ams}            {}{American Mathematical Society}{Providence}
\DefinePublisher{academiclondon} {}{Academic Press}               {London}
\DefinePublisher{academicnewyork}{}{Academic Press}               {New~York}
\DefinePublisher{dekker}         {}{Dekker}                       {New~York}
\DefinePublisher{springer}       {}{Springer}                     {Berlin}
\DefinePublisher{kluwer}         {}{Kluwer}                       {Dordrecht}
\DefinePublisher{cambridge}      {}{Cambridge University Press}   {}
\DefinePublisher{scottforesman}  {}{Scott, Foresman \& Co.}       {Chicago}
\DefinePublisher{cremonese}      {}{Cremonese}                    {Roma}
%
\bib{01}{article}{
 title={Multipliers of $p$-integrable functions}*{language={english}},
 journal={Bulletin of the American Mathematical Society}*{language={english}},
 volume={70},
 date={1964},
 pages={666--669},
 doi={10.1090/S0002-9904-Sotirios1964-11153-8},
}\pages{666}{669}
%
\bib{02}{article}{
 title={Translation invariant operators in $L^p$}*{language={english}},
 journal={Duke Mathematical Journal}*{language={english}},
 volume={32},
 date={1965},
 pages={495--501},
}\pages{495}{501}
%
\bib{03}{article}{
 title={On the subspace of $L^p$ invariant under multiplication of transform by bounded continuous functions}*{language={english}},
 journal={Rendiconti del Seminario Matematico della Università di Padova}*{language={italian}},
 volume={35},
 date={1965},
 number={1},
 pages={176--189},
}\pages{176}{189}
%
\bib{04}{article}{
 author={curtis}*{language={english}},
 title={Factorization theorems for Banach algebras}*{language={english}},
 conference={
  title={Function algebras}*{language={english}},
  address={Tulane University}*{language={english}},
  date={1965},
 },
 book={
  publisher={scottforesman}*{language={english}},
  date={1966},
 },
 pages={169--185},
}\pages{169}{185}
%
\bib{05}{article}{
 author={rider}*{language={english}},
 title={A theorem of Littlewood and lacunary series for compact groups}*{language={english}},
 journal={Pacific Journal of Mathematics}*{language={english}},
 volume={16},
 date={1966},
 number={3},
 pages={505--514},
}\pages{505}{514}
%
\bib{06}{article}{
 author={gaudry}*{language={english}},
 title={Density and representation theorems for multipliers of type $(p,q)$}*{language={english}},
 journal={Journal of the Australian Mathematical Society}*{language={english}},
 volume={7},
 date={1967},
 pages={1--6},
}\pages{1}{6}
%
\bib{07}{article}{
 author={rider}*{language={english}},
 title={A theorem on random Fourier series on noncommutative groups}*{language={english}},
 journal={Pacific Journal of Mathematics}*{language={english}},
 volume={21},
 date={1967},
 number={3},
 pages={487--492},
}\pages{487}{492}
%
\bib{08-Ex09}{article}{
 title={Appartenenza a $L^p$ delle serie di Fourier aleatorie su gruppi non commutativi}*{language={italian}},
 journal={Rendiconti del Seminario Matematico della Università di Padova}*{language={italian}},
 volume={39},
 date={1967},
 pages={330--348},
}\pages{330}{348}
%
\bib{09-Ex08}{article}{
 title={Un'osservazione sulle serie di Fourier lacunari per gruppi non commutativi}*{language={italian}},
 journal={Bollettino dell'Unione Matematica Italiana, serie~3,}*{language={italian}},
 volume={22},
 date={1967},
 number={4},
 pages={497--504},
}\pages{497}{504}
%
\bib{10}{article}{
 title={Una maggiorazione a priori per equazioni ellittiche in due variabili}*{language={italian}},
 journal={Bollettino dell'Unione Matematica Italiana, serie~4,}*{language={italian}},
 volume={1},
 date={1968},
 pages={107--121},
}\pages{107}{121}
%
\bib{11}{article}{
 title={Bounded and continuous random Fourier series on non-commutative groups}*{language={english}},
 journal={Proceedings of the American Mathematical Society}*{language={english}},
 volume={22},
 date={1969},
 number={3},
 pages={573--578},
 doi={10.2307/2037434},
}\pages{573}{578}
%
\bib{12}{article}{
 author={gaudry}*{language={english}},
 title={Multipliers and sets of uniqueness of $L^p$}*{language={english}},
 journal={Michigan Mathematical Journal}*{language={english}},
 volume={17},
 date={1970},
 pages={179--191},
}\pages{179}{191}
%
\bib{13}{article}{
 title={An example in the theory of lacunary Fourier series}*{language={english}},
 journal={Bollettino dell'Unione Matematica Italiana, serie~4,}*{language={italian}},
 volume={3},
 date={1970},
 pages={375--378},
}\pages{375}{378}
%
\bib{14}{article}{
 author={gaudry}*{language={english}},
 title={Extensions of multipliers}*{language={english}},
 journal={Bollettino dell'Unione Matematica Italiana, serie~4,}*{language={italian}},
 volume={3},
 date={1970},
 pages={1003--1014},
}\pages{1003}{1014}
%
\bib{15}{article}{
 author={franklin}*{language={english}},
 title={Multipliers of distributive lattices}*{language={english}},
 journal={Indian Journal of Mathematics}*{language={english}},
 volume={12},
 date={1970},
 number={3},
 pages={153--161},
}\pages{153}{161}
%
\bib{16}{article}{
 author={gaudry}*{language={english}},
 title={Multipliers of $L^p$ which vanish at infinity}*{language={english}},
 journal={Journal of Functional Analysis}*{language={english}},
 volume={7},
 date={1971},
 pages={475--486},
 doi={10.1016/0022-1236(71)90029-2},
}\pages{475}{486}
%
\bib{17}{article}{
 title={Random Fourier series on compact groups}*{language={english}},
 conference={
  title={Theory of group representations and Fourier analysis, II Ciclo}*{language={english}},
  address={Montecatini Terme}*{language={italian}},
  date={1970},
 },
 book={
  series={Centro Internazionale Matematico Estivo}*{language={italian}},
  publisher={cremonese}*{language={italian}},
  date={1971},
 },
 pages={1--63},
}\pages{1}{62}
%
\bib{18}{article}{
 title={Recensione dell'opera
 \textit{Abstract harmonic analysis}, voll. I e II, di Edwin Hewitt e Kenneth~A. Ross}*{language={english}},
 journal={Bulletin of the American Mathematical Society}*{language={english}},
 volume={78},
 date={1972},
 number={2},
 pages={172--178},
}\pages{172}{178}
%
\bib{19}{article}{
 author={price}*{language={english}},
 title={Applications of random Fourier series over compact groups to Fourier multipliers}*{language={english}},
 journal={Pacific Journal of Mathematics}*{language={english}},
 volume={43},
 date={1972},
 number={2},
 pages={531--541},
}\pages{531}{541}
%
\bib{20}{article}{
 author={price}*{language={english}},
 title={Rudin-Shapiro sequences on compact groups}*{language={english}},
 journal={Bulletin of the Australian Mathematical Society}*{language={english}},
 volume={8},
 date={1973},
 pages={241--245},
 doi={10.1017/S0004972700042490},
}\pages{241}{245}
%
\bib{21}{article}{
 author={picardello}*{language={italian}},
 title={Multiplicateurs de $A(G)$ qui ne sont pas dans $B(G)$}*{language={french}},
 journal={Comptes Rendus Hebdomadaires des Séances de l'Académie des Sciences, serie~A: Sciences Mathématiques}*{language={french}},
 volume={277},
 date={1973},
 pages={117--119},
}\pages{117}{119}
%
\bib{22}{article}{
 author={cecchini}*{language={italian}},
 title={Projections of uniqueness for $L^p(G)$}*{language={english}},
 journal={Pacific Journal of Mathematics}*{language={english}},
 volume={51},
 date={1974},
 number={1},
 pages={37--47},
}\pages{37}{47}
%
\bib{23}{article}{
 title={Multipliers vanishing at infinity for certain compact groups}*{language={english}},
 journal={Proceedings of the American Mathematical Society}*{language={english}},
 volume={45},
 date={1974},
 number={2},
 pages={199--203},
 doi={10.2307/2040062},
}\pages{199}{203}
%
\bib{24}{article}{
 title={Positive definite functions which vanish at infinity}*{language={english}},
 journal={Pacific Journal of Mathematics}*{language={english}},
 volume={69},
 date={1977},
 number={2},
 pages={355--363},
}\pages{355}{363}
%
\bib{25}{article}{
 title={On the action of unitary groups on a Hilbert space}*{language={english}},
 conference={
  title={Analisi armonica e spazi di funzioni su gruppi localmente compatti}*{language={italian}},
  address={Roma}*{language={italian}},
  date={1976},
 },
 book={
  series={Istituto Nazionale di Alta Matematica, Symposia Mathematica}*{language={italian}},
  volume={22},
  publisher={academiclondon}*{language={english}},
  date={1977},
 },
 pages={315--319},
}\pages{315}{319}
%
\bib{26}{article}{
 title={Insiemi lacunari nei gruppi non commutativi}*{language={italian}},
 journal={Rendiconti del Seminario Matematico e Fisico di Milano}*{language={italian}},
 volume={47},
 date={1977},
 pages={45--59 (1979)},
 doi={10.1007/BF02925741},
}\pages{45}{59}
%
\bib{27}{article}{
 author={picardello}*{language={italian}},
 title={Functions that operate on the algebra $B_0(G)$}*{language={english}},
 journal={Pacific Journal of Mathematics}*{language={english}},
 volume={74},
 date={1978},
 number={1},
 pages={57--61},
}\pages{57}{61}
%
\bib{28}{article}{
 author={mauceri}*{language={italian}},
 title={A characterization of nonatomic Hilbert algebras}*{language={english}},
 journal={Proceedings of the American Mathematical Society}*{language={english}},
 volume={72},
 date={1978},
 number={3},
 pages={468--472},
 doi={10.2307/2042453},
}\pages{468}{472}
%
\bib{29}{article}{
 title={Construction of positive definite functions}*{language={english}},
 journal={Rendiconti del Seminario Matematico della Università di Padova}*{language={italian}},
 volume={60},
 date={1978},
 pages={43--53 (1979)},
}\pages{43}{53}
%
\bib{30}{article}{
 title={A remark on multipliers of the Fourier algebra of the free group}*{language={english}},
 journal={Bollettino dell'Unione Matematica Italiana, serie~5, sezione~A,}*{language={italian}},
 volume={16},
 date={1979},
 number={3},
 pages={577--581},
}\pages{577}{581}
%
\bib{31}{article}{
 title={Singular positive definite functions}*{language={english}},
 conference={
  title={Harmonic analysis}*{language={english}},
  address={Iraklion},
  date={1978},
 },
 book={
  series={Lecture Notes in Mathematics}*{language={english}},
  volume={781},
  publisher={springer}*{language={german}},
  date={1980},
 },
 pages={22--29},
}\pages{22}{29}
%
\bib{32}{article}{
 author={demichele}*{language={italian}},
 title={Positive definite functions on free groups}*{language={english}},
 journal={American Journal of Mathematics}*{language={english}},
 volume={102},
 date={1980},
 number={3},
 pages={503--509},
 doi={10.2307/2374112},
}\pages{503}{509}
%
\bib{33}{article}{
 author={picardello}*{language={italian}},
 title={Spherical functions and harmonic analysis on free groups}*{language={english}},
 journal={Journal of Functional Analysis}*{language={english}},
 volume={47},
 date={1982},
 number={3},
 pages={281--304},
 doi={10.1016/0022-1236(82)90108-2},
}\pages{281}{304}
%
\bib{34}{book}{
 author={picardello}*{language={italian}},
 title={Harmonic analysis on free groups}*{language={english}},
 series={Lecture Notes in Pure and Applied Mathematics}*{language={english}},
 volume={87},
 publisher={dekker}*{language={english}},
 date={1983},
}\pages{1}{268}
%
\bib{35}{article}{
 author={picardello},
 title={Restriction of spherical representations of $\operatorname{PGL}_2(\mathbf{Q}_p)$ to a discrete subgroup}*{language={english}},
 journal={Proceedings of the American Mathematical Society}*{language={english}},
 volume={91},
 date={1984},
 number={3},
 pages={405--408},
 doi={10.2307/2045312},
}\pages{405}{408}
%
\bib{36}{article}{
 title={Analisi armonica su strutture discrete}*{language={italian}},
 journal={Bollettino dell'Unione Matematica Italiana, serie~6, sezione~A,}*{language={italian}},
 volume={3},
 date={1984},
 number={3},
 pages={313--334},
}\pages{313}{334}
%
\bib{37}{article}{
 author={steger},
 title={Harmonic analysis on trees}*{language={english}},
 conference={
  title={Harmonic analysis, symmetric spaces and probability theory}*{language={english}},
  address={Cortona}*{language={italian}},
  date={1984},
 },
 book={
  series={Symposia Mathematica}*{language={italian}},
  volume={29},
  publisher={academicnewyork}*{language={english}},
  date={1987},
 },
 pages={163--182},
}\pages{163}{182}
%
\bib{38}{article}{
 author={cohen}*{language={english}},
 title={Idempotents in the reduced $C^*$-algebra of a free group}*{language={english}},
 journal={Proceedings of the American Mathematical Society}*{language={english}},
 volume={103},
 date={1988},
 number={3},
 pages={779--782},
 doi={10.2307/2046852},
}\pages{779}{782}
%
\bib{39}{book}{
 author={nebbia}*{language={italian}},
 title={Harmonic analysis and representation theory for groups acting on homogeneous trees}*{language={english}},
 series={London Mathematical Society Lecture Note Series}*{language={english}},
 volume={162},
 publisher={cambridge}*{language={english}},
 date={1991},
 doi={10.1017/CBO9780511662324},
}\pages{1}{156}
%
\bib{40}{article}{
 author={steger}*{language={english}},
 title={Harmonic analysis for anisotropic random walks on homogeneous trees}*{language={english}},
 journal={Memoirs of the American Mathematical Society}*{language={english}},
 volume={110},
 date={1994},
 number={531},
 doi={10.1090/memo/0531},
}\pages{1}{74}
%
\bib{41}{article}{
 title={Diffusion on compact ultrametric spaces}*{language={english}},
 conference={
  title={Noncompact Lie groups and some of their applications}*{language={english}},
  address={San~Antonio}*{language={spanish}},
  date={1993},
 },
 book={
  series={NATO Advanced Science Institutes Series C: Mathematical and Physical Sciences}*{language={english}},
  volume={429},
  publisher={kluwer}*{language={dutch}},
  date={1994},
 },
 pages={157--167},
}\pages{157}{167}
%
\bib{42}{article}{
 title={Local fields and trees}*{language={english}},
 conference={
  title={Harmonic functions on trees and buildings}*{language={english}},
  address={New York}*{language={english}},
  date={1995},
 },
 book={
  series={Contemporary Mathematics}*{language={english}},
  volume={206},
  publisher={ams}*{language={english}},
  date={1997},
 },
 pages={3--16},
 doi={10.1090/conm/206/02684},
}\pages{3}{16}
%
\bib{43}{article}{
 author={baldi},
 author={casadiotarabusi}*{language={italian}},
 title={Stable laws arising from hitting distributions of processes on homogeneous trees and the hyperbolic half-plane}*{language={english}},
 journal={Pacific Journal of Mathematics}*{language={english}},
 volume={197},
 date={2001},
 number={2},
 pages={257--273},
 doi={10.2140/pjm.2001.197.257},
}\pages{257}{273}
%
\bib{44}{article}{
 title={An application of Gelfand pairs to a problem of diffusion in compact ultrametric spaces}*{language={english}},
 conference={
  title={Topics in probability and Lie groups: boundary theory}*{language={english}},
  address={Centre de Recherches Mathématiques},
 },
 book={
  series={CRM Proceedings and Lecture Notes}*{language={english}},
  volume={28},
  publisher={ams}*{language={english}},
  date={2001},
 },
 pages={51--67},
 doi={10.1090/crmp/028/03},
}\pages{51}{67}
%
\bib{45}{article}{
 author={baldi}*{language={italian}},
 author={casadiotarabusi}*{language={italian}},
 author={yor},
 title={Non-symmetric hitting distributions on the hyperbolic half-plane and subordinated perpetuities}*{language={english}},
 journal={Revista Matemática Iberoamericana}*{language={spanish}},
 volume={17},
 date={2001},
 number={3},
 pages={587--605},
 doi={10.4171/RMI/305},
}\pages{587}{605}
%
\bib{46}{article}{
 author={delmuto},
 title={Diffusion on locally compact ultrametric spaces}*{language={english}},
 journal={Expositiones Mathematicae}*{language={latin}},
 volume={22},
 date={2004},
 number={3},
 pages={197--211},
 doi={10.1016/S0723-0869(04)80005-7},
}\pages{197}{211}
%
\bib{47}{article}{
 author={delmuto},
 title={Anisotropic diffusion on totally disconnected abelian groups}*{language={english}},
 journal={Pacific Journal of Mathematics}*{language={english}},
 volume={225},
 date={2006},
 number={2},
 pages={221--229},
 doi={10.2140/pjm.2006.225.221},
}\pages{221}{229}
%
\bib{48}{article}{
 author={casadiotarabusi},
 title={Drifted Laplace operators on homogeneous trees}*{language={english}},
 journal={Proceedings of the American Mathematical Society}*{language={english}},
 volume={135},
 date={2007},
 number={7},
 pages={2165--2175},
 doi={10.1090/S0002-9939-07-08811-9},
}\pages{2165}{2175}
%
\bib{49}{article}{
 author={guerriero}*{language={italian}},
 author={leone}*{language={italian}},
 author={mignoli}*{language={italian}},
 author={rogora}*{language={italian}},
 title={Decomposition of variance in terms of conditional means}*{language={english}},
 journal={Statistica (Bologna)}*{language={italian}},
 volume={67},
 date={2007},
 number={2},
 pages={191--201 (2008)},
}\pages{191}{201}
%
\bib{50}{article}{
 author={casadiotarabusi}*{language={italian}},
 title={Poisson kernels of drifted Laplace operators on trees and on the half-plane}*{language={english}},
 journal={Colloquium Mathematicum}*{language={latin}},
 volume={118},
 date={2010},
 number={1},
 pages={147--159},
 doi={10.4064/cm118-1-7},
}\pages{147}{159}
%
\end{biblist}
%
\end{bibdiv}
\end{document}